\documentclass[review,10pt,authoryear]{elsarticle}
\usepackage[margin=2.5cm]{geometry}
\usepackage{graphicx} 
\usepackage[hidelinks]{hyperref}
\usepackage{rotating}
\usepackage{amssymb, amsmath, amsfonts}
\usepackage{mathrsfs}
\usepackage{caption}
\usepackage{float}
\usepackage{comment}
\usepackage{lscape}
\usepackage{cancel}
\usepackage{autobreak}
\usepackage{algorithmicx}
\usepackage{algcompatible}
\usepackage{framed} 
\usepackage{longtable} 
\usepackage{xcolor}
\usepackage{relsize}
\usepackage{subfig}
\usepackage[linesnumbered,ruled,vlined,boxed,longend]{algorithm2e}
\SetStartEndCondition{ }{ }{ }
\SetKwFor{ForEach}{for each}{do}{end for}
\SetArgSty{textnormal}

\usepackage{algpseudocode}
\usepackage{amsmath, amssymb}
\usepackage{algpseudocode}
\usepackage{tcolorbox}
\usepackage{xcolor}
\usepackage{booktabs}
\usepackage{amsthm}
\usepackage{comment}
\usepackage{eurosym}
\usepackage{float}
\usepackage{booktabs}
\usepackage{multirow}
\usepackage[normalem]{ulem}
\usepackage{url}
\usepackage{lineno}

\begin{document}
\begin{frontmatter}
\title{Roadside Units Placement for Total Travel Time Minimization in I2X-Enabled Road Networks Based on a Dynamic Traffic Flow Model
}
\author[a]{Emiliano Cristiani} 
\author[a,b]{Francesca L. Ignoto}
\author[c,d]{Anna Laura Pala\corref{cor1}}
\author[d]{Giuseppe Stecca}
\address[a]{Istituto per le Applicazioni del Calcolo, Consiglio Nazionale delle Ricerche, Via dei Taurini 19, 00185 Rome, Italy}
\address[b]{Dipartimento di Scienze di Base e Applicate per l'Ingegneria, Sapienza Università di Roma, Via Scarpa 16, 00161 Rome, Italy}
\address[c]{Dipartimento di Ingegneria Informatica, Automatica e Gestionale, Sapienza Università di Roma, Via Ariosto 25, Rome, 00185, Italy}
\address[d]{Istituto di Analisi dei Sistemi ed Informatica ``A. Ruberti'', Consiglio Nazionale delle Ricerche, Via dei Taurini 19, 00185 Rome, Italy}

\begin{abstract}
This work addresses the problem of locating Road Side Units (RSUs) in vehicular networks with the aim of improving overall traffic efficiency. We propose a heuristic approach for determining the placement of RSUs on a road network, while vehicle dynamics are modeled through a microscopic
follow-the-leader interaction in which drivers continuously adapt their speed according to the distance from the vehicle ahead. The RSUs share the collected information among themselves to cooperatively estimate the network state,  which is then communicated to vehicles, allowing them to compute their fastest paths. We also assume that not all drivers comply with the provided recommendations, meaning that only a fraction of users follow the suggested routes. Since the aim is to improve network performance while minimizing the number of deployed units, the methodology adopts a nested optimization scheme: an outer module explores the optimal number $k$ of RSUs, while an inner module identifies the best placement configuration for each $k$. To evaluate the effectiveness and the scalability of the proposed method, we tested two networks with different topologies and sizes. The results show a decrease in Total Travel Time ($TTT$) compared to the baseline scenario (i.e., without RSUs) of up to 36\% on the smaller network and up to 20\% on the larger one. Stochastic vehicle origin-destination pairs have also been considered.
\end{abstract}


\begin{keyword}
RSUs placement;
Infrastructure-to-Everything (I2X) communication;
Traffic Optimization;
Second-order dynamic traffic models;
Driver compliance;
Road networks.
\end{keyword}
\end{frontmatter}
\section{Introduction}
\paragraph{Context} In recent years, widespread urbanization has significantly increased urban traffic congestion, introducing critical challenges not only for environmental sustainability but also for public safety due to a corresponding increase in traffic accidents \citep{WHO2018RoadSafety}. In this context, the deployment of connected and Intelligent Transportation Systems (ITS) has emerged as a promising research area to improve the efficiency and safety of transportation networks, enabling data-driven and adaptive traffic management \citep{ 5430544, silva2017survey, elassy2024intelligent}. 

Within ITS, route planning and dynamic traffic guidance play a key role in efficiently redistributing traffic flows, alleviating congestion, and reducing network travel times \citep{Stecca26, PALA2026102002}. To enable accurate traffic‑flow prediction, Vehicular Ad‑Hoc Networks (VANETs) and their evolution in Internet of Vehicles (IoV) are essential, as they provide real‑time traffic information in highly dynamic and mobile network environments \citep{akhtar2014vehicle}. VANET is a wireless network that is based primarily on vehicle-to-vehicle (V2V) and vehicle-to-infrastructure (V2I) communication. Specifically, V2V ensures message transmission when two or more vehicles are in the same transmission range, and even beyond that range through multi-hopping \citep{heidari2014survey}. V2I communication allows real-time traffic information to be collected, processed, and redistributed through roadside devices such as Road Side Units (RSUs) \citep{hartenstein2009vanet, sommer2014vehicular}. Information between road side sensors and vehicles can be exchanged either using short-range communications or cellular communications. 

These systems have the potential to improve traffic efficiency by providing drivers with routing recommendations based on the current state of the network. However, their effectiveness depends on several factors, including the topology of the road network, the characteristics of specific road segments (e.g., accident-prone areas), and the spatial and temporal distribution of vehicles, as highlighted in the survey by \cite{guerna2022roadside}. Despite their robustness, these devices cannot be deployed indefinitely due to their high installation and maintenance costs \citep{nikookaran2017combining}. This limitation suggests the critical importance of the RSUs placement optimization problem, which aims to maximize network performance and communication coverage while minimizing costs.
%
Several studies in literature have addressed RSUs placement optimization, aiming to improve a wide range of performance indicators (e.g., user service levels, network coverage and communication latency). 
However, most approaches rely on static traffic models and fail to capture realistic traffic dynamics, as they do not account for variations in vehicle speed and acceleration. In this study, vehicle dynamics are instead described through a microscopic follow-the-leader interaction, in which drivers continuously adapt their speed according to the distance from the vehicle ahead. This model reproduces car-to-car interactions and captures complex traffic phenomena such as stop-and-go waves, allowing us to investigate the impact of infrastructure-to-vehicle (I2V) communications on realistic traffic dynamics. 

In our framework, I2V communications occur whenever a vehicle passes next to a RSU: the RSU measures the current vehicle velocity and, in exchange, provides to vehicle necessary travel time information for routing optimization. 
We also assume that RSUs are interconnected through infrastructure-to-infrastructure (I2I) communication, which allows them to aggregate the collected data and estimate current travel times on the network. 
Further, we assume that drivers are not fully compliant with the suggested routes. In practice, only a fraction of users follows the routing recommendations provided by the infrastructure, while the remaining drivers keep their original route choice. Modeling partial compliance is important for realistically assessing the impact of information systems on traffic dynamics \citep{bonsall1991, ben2001route}. 

\medskip

\paragraph{Main contribution} We propose a novel methodology for optimizing RSUs placement, designed as an iterative feedback loop in which candidate locations for RSUs are evaluated through second-order (i.e., acceleration-based) dynamic traffic simulations.
While most of the literature aims to determine the placement of RSUs for maximizing coverage and network communication indicators, our work focuses on determining the optimal number of RSUs and their location in order to minimize the Total Travel Time ($TTT$). 
To address the problem, the proposed location algorithm is based on a specialized Iterated Local Search metaheuristic, in which candidate configurations are generated and iteratively improved using a dual-criterion approach: the path road coverage, which measures the number of different routes that traverse a link, and the link origin-destination (OD) coverage, which quantifies the number of distinct downstream OD flows intercepted by a link. The effectiveness of each configuration is evaluated through traffic simulations, using $TTT$ as a network-level metric of traffic efficiency.
A second metric, i.e.\ the total \emph{fuel consumption} estimation, is also considered for assessing the environmental impact of the proposed approach. 

\medskip 

The remainder of the paper is as follows: Section~\ref{sec:literature} discusses relevant literature; Section~\ref{sec:methodology} defines the micro-simulation model and the RSUs optimization framework;  Section~\ref{sec:results} describes the considered road networks and computational setup, and analyzes the numerical results.
Section~\ref{sec:conclusion} concludes the paper.

\section{Literature Review} \label{sec:literature}
Within the architecture of a VANET, vehicles are equipped with On-board Processing Units (OPUs), while Roadside Units (RSUs) are strategically deployed along highways and at road intersections to provide traffic directives \citep{guerna2022roadside}. 
Through these nodes, vehicles dynamically exchange traffic information (e.g., position, speed, local traffic density) with nearby vehicles using V2V communications, or with infrastructure elements such as RSUs via Vehicle-to-RSU (V2R or V2I) communications.

Conventional routing and navigation systems (e.g., Google Maps, Waze) typically rely on centralized cloud servers to evaluate optimal travel paths. However, in urban settings these centralized configurations face operational challenges, including severe scalability constraints, end-to-end communication latencies, and vulnerabilities associated with single points of failure. 
To address these architectural limitations and satisfy the requirements of autonomous driving, e.g., low-latency communication \citep{ge2019ultra}, current research shifts toward decentralized edge computing systems where data is processed directly on the road network by local RSUs.
As robust nodes placed right next to the drivers, RSUs can operate independently even during total network outages and can improve the connectivity of the network by dramatically reduce communication delays \citep{ni2018joint}. 
Therefore, strategic RSUs placement is crucial to improve the performance of vehicular networks, whether in static \citep{guerna2019gica} or dynamic deployment settings \citep{cai2020combined}. RSUs placement naturally falls within the broader class of facility location problems, a family of optimization models extensively studied since the 1960s in numerous variants \citep{owen1998strategic}. When considering RSUs deployment, the objectives usually focus either on improving communication performance or on reducing the infrastructure costs involved in construction, installation, and maintenance. Accordingly, several approaches have been presented in the literature, including data-driven techniques, single-objective optimization, and multi-objective optimization. 

A large portion of the literature on RSUs deployment relies on single-objective optimization, using either heuristic methods or mathematical programming. For instance, several studies focus exclusively on maximizing network coverage \citep {cai2020trajectory} and information dissemination \citep{eftekhari2015binary}, while others on maximizing service profit \citep{ni2018joint} or minimizing costs \citep{guerna2021ac}. Within this context, \cite{trullols2010planning} formulated RSUs placement as a maximum‑coverage optimization problem and developed heuristic algorithms to solve it efficiently under realistic urban mobility conditions. \cite{aslam2012optimal} developed a so-called balloon optimization heuristic that models RSUs coverage as expanding regions in a 2‑D space, achieving optimal or near‑optimal RSUs placements compared with exhaustive search. \cite{cavalcante2012roadside} modeled the problem as a Maximum Coverage with Time Threshold Problem (MCTTP) and used a genetic algorithm to solve it.  \cite{ghorai2018constrained} developed a constrained Delaunay Triangulation‑based method for RSUs placement that incorporates obstacles and jointly optimizes node locations and transmission ranges in complex urban environments, while \cite{10.1145/2507248.2507250} introduced a geometry-driven coverage approach to address the RSUs deployment problem in urban settings. \cite{eftekhari2015binary} introduced a binary programming (BP) model for maximizing information dissemination to vehicles. 
More recently, \cite{qi2020traffic} modeled data‑delivery optimization as a 0-1 Non-Linear Integer Programming (NLIP) and designed a two‑stage heuristic for traffic‑differentiated clustering and routing that minimizes end‑to‑end delay while controlling cellular bandwidth usage. 

Recent literature has shifted toward multi-objective optimization models. For example, \cite{guerna2019gica} formulated RSUs placement as a multi‑objective optimization model and proposed a genetic intersection‑coverage heuristic that prioritizes high‑impact intersections to maximize connectivity while reducing interference and deployment cost. 
\cite{liang2023deploying} introduced a  multi‑objective model maximizing served tasks and minimizing task‑weighted delay, reformulated it as a p‑median problem, and solved it through a multi‑objective Lagrangian‑relaxation heuristic. Extending this paradigm, \cite{Guo2025} proposed two multi‑objective evolutionary algorithms that handle obstacles and large decision spaces via multi‑population search and adaptive exploration, with RSU‑density calibration and an iterative best‑response game for data offloading.

Despite the extensive literature on RSUs deployment, only a few studies explicitly address placement strategies aimed at reducing congestion or $TTT$. Several works address re‑routing and route‑guidance mechanisms in VANETs, but they typically operate over a fixed RSUs infrastructure and do not optimize RSUs placement itself. For example, \cite{dong2026ccrgs} designed a cluster‑based collaborative route‑guidance strategy that balances road congestion and RSU computational load through a multi‑factor road‑weight model, a bi‑objective A* routing algorithm, vehicle clustering, and periodic dynamic re‑routing, assuming a fixed RSUs deployment. \cite{tay2025intelligent} implemented a centralized re‑routing system where fixed RSUs collect traffic data for the TMC, which uses an SVM model to predict near‑future traffic conditions, evaluates congestion through a custom KPI, selects vehicles for diversion, and computes alternative routes via a k‑shortest‑path algorithm to reduce $TTT$. \cite{guo2018real} designed a VANET‑assisted real‑time path‑planning system that relies on RSU‑based information sharing to estimate travel times and dynamically reroute vehicles, using a distributed V2V/V2R/R2R communication architecture and a TTE‑driven path‑selection algorithm to avoid congestion. \cite{talusan2020route} developed a decentralized route‑planning service that runs entirely on RSUs, using a distributed edge‑computing middleware (RSU‑Edge) and a task‑allocation algorithm that balances accuracy and response time by reallocating routing queries among neighboring RSUs. Finally, \cite{guo2020dynamic} proposed an analytic optimization model that minimizes $TTT$ and traffic‑load imbalance, and solved it through a Dynamic Interior Point Method (DIPM) that reroutes drivers via a central server using real‑time traffic information collected through cellular networks or VANETs with RSUs.

\medskip 

Several gaps emerge from the literature analysis. First, the vast majority of RSUs placement studies adopt coverage-oriented or communication-centric objectives, while our proposal addresses optimal location of RSUs with the aim of $TTT$ minimization. Second, nearly all approaches employ static or aggregated traffic representations and not second-order (i.e.\ acceleration-based) dynamic models. Third, no existing RSUs placement framework integrates the microscopic traffic model within a closed-loop specific optimization algorithm, considering path-road coverage and link origin-destination coverage criteria. Fourth, with respect to the literature, our approach also considers the environmental impact of RSUs placements by accounting for total fuel consumption.

\section{Methodology} \label{sec:methodology}
The proposed approach iteratively combines microscopic dynamic traffic simulation with a multi-criteria optimization framework to determine the optimal RSU placement. The framework is designed as an iterative feedback loop where candidate configurations generated by the optimization module are evaluated through traffic simulation to quantify their impact on network performance. In particular, for each candidate configuration, vehicle traffic is described through a microscopic follow-the-leader interaction, in which drivers adapt their speed according to the distance from the vehicle ahead (second-order dynamics)). The formulation follows the model used in \cite{briani2021macroscopic}, in turn inspired by the model proposed by \cite{zhao2017unified}. 

Each RSU is equipped with traffic monitoring sensors that measure the speed of vehicles as they pass beside the RSU. These measurements are used to estimate the travel time on the monitored road, providing traffic information for dynamic route guidance.

The proposed methodology relies on two communication technologies: 
\begin{itemize}
    \item \textbf{I2I communication:} RSUs continuously exchange the travel time information they collect. As a result, each RSU progressively builds a global knowledge of the traffic conditions over the network;
    \item \textbf{I2V communication:} whenever a vehicle passes an RSU, the infrastructure transmits the latest travel time estimates available for all monitored roads. 
    Compliant vehicles exploit this information to dynamically recompute their shortest path towards the destination.
\end{itemize}

 To ensure consistency throughout the rest of this paper, we define the following notation. The road network is modeled as a directed graph $G = (V,L)$, where $V$ is the set of nodes (i.e., road intersections) and $L$ denotes the set of links (i.e., road segments). 
 A RSU configuration $\mathbf{x} \in \{0,1\}^{|L|}$ is a binary vector with entries equal to 1 if the corresponding links are equipped with RSUs and 0 otherwise. The performance of a given RSU configuration $\mathbf{x}$ with $k$ units is quantified by its impact on the network's $TTT$, denoted by $TTT_{\mathbf{x}}$. 
In particular, our aim is to determine the RSU configuration $\mathbf{x}$
that maximizes the objective $\delta_{TTT}(\mathbf{x})$, representing the relative reduction of the network’s $TTT$ compared to the baseline scenario $TTT_0$ (i.e., the $TTT$ in the absence of any RSUs):
\begin{equation}
\delta_{TTT}(\mathbf{x}) = \frac{TTT_0 - TTT_{\mathbf{x}}}{TTT_0}.
\end{equation}
We further define the global performance bounds on the network by considering all possible RSU configurations:
\begin{equation}
TTT_{\min} = \min_{\mathbf{x}} TTT_{\mathbf{x}},
\end{equation}
\begin{equation}
TTT_{\max} = \max_{\mathbf{x}} TTT_{\mathbf{x}}.
\end{equation}
Accordingly, we also define the global reduction metric $\Delta_{TTT}$, which quantifies the maximum improvement in $TTT$ compared to the baseline scenario:
\begin{equation}
\Delta_{TTT} = \frac{TTT_0 - TTT_{\min}}{TTT_0}.
\end{equation}
Finally, to quantify the global range of network performance, we introduce the parameter $\Gamma$ as follows
\begin{equation}
\Gamma_{TTT} = \frac{TTT_{\max} - TTT_{\min}}{TTT_{\max}}.
\end{equation}
In the following, the two modules of the proposed methodology are presented in detail.

\subsection{Microscopic traffic model} \label{sec:traffic_simulation}
We consider a set of $N_c$ vehicles traveling over a road network. Each vehicle is assigned a fixed origin-destination pair together with a compliance status, which specifies whether the driver exploits the routing recommendations provided by the infrastructure. A compliance rate $\gamma$ is fixed to represent the fraction of compliant vehicles, and the corresponding vehicles are randomly selected at the beginning of the simulation.

\subsubsection{Vehicle dynamics}\label{sec:dynamics}
Let $X_c(t)$ and $V_c(t)$ denote the position and velocity of vehicle $c$ at time $t$, respectively, for $c=1,\ldots,N_c$. The evolution of the system is governed by the following system of ordinary differential equations
\begin{equation}
\begin{array}{l}
\dot{X}_c(t)=V_c(t), \\
\dot{V}_c(t)=A\left(X_c,X_{c'},V_c,V_{c'}\right),
\end{array}
\end{equation}
where the acceleration function $A$ depends on the state of the leading vehicle $c'$, i.e., the vehicle immediately preceding vehicle $c$ \emph{along its route} at time $t$. Denoting by $\delta_c(t)$ the distance between the vehicle $c$ and the vehicle $c'$ at time $t$, the acceleration is modeled as a relaxation toward an equilibrium velocity
$v_{eq}(\delta_c)$, given by
\begin{equation}\label{eq:acc}
A_c:=A\left(X_c,X_{c'},V_c,V_{c'}\right)=
\begin{cases}
\dfrac{1}{\tau_{\mathrm{acc}}}
\left(v_{eq}(\delta_c)-V_c\right),
& \text{if } v_{eq}(\delta_c)\geq V_c,
\\[0.3cm]
\dfrac{1}{\tau_{\mathrm{dec}}}
\left(v_{eq}(\delta_c)-V_c\right),
& \text{if } v_{eq}(\delta_c)<V_c,
\end{cases}
\end{equation}
where the parameters $\tau_{\mathrm{acc}}$ and $\tau_{\mathrm{dec}}$ represent
the acceleration and deceleration reactivity times,
respectively. Using different values for the two parameters allows the model
to reproduce the asymmetric response of drivers, who generally brake more
rapidly than they accelerate.

The equilibrium velocity is defined as a piecewise-linear function of the
headway:
\begin{equation}\label{eq:veq}
v_{eq}(\delta)=
\begin{cases}
0,
& \delta \leq \delta_{\mathrm{close}},
\\
V_{\max}
\dfrac{\delta-\delta_{\mathrm{close}}}
{\delta_{\mathrm{far}}-\delta_{\mathrm{close}}},
& \delta_{\mathrm{close}}<\delta<\delta_{\mathrm{far}},
\\[0.1cm]
V_{\max},
& \delta \geq \delta_{\mathrm{far}},
\end{cases}
\end{equation}
where $\delta_{\mathrm{close}}$ is the minimum safe distance,
$\delta_{\mathrm{far}}$ is the distance beyond which drivers can travel at
their desired speed, and $V_{\max}$ denotes the maximum admissible velocity. According to \eqref{eq:veq}, the desired velocity increases with the available headway, ranging from zero in highly congested conditions to the free-flow speed when sufficient space is available ahead. This interaction mechanism provides a realistic representation of traffic dynamics and captures phenomena such as queue formation and stop-and-go waves  \citep{briani2021macroscopic}.

\subsubsection{Path choice}
Before the simulation starts, each vehicle computes the route connecting its assigned origin and destination through its OPU. Following the destination-aware routing framework proposed in \cite{cristiani2026microscopic}, the route is determined as the shortest path assuming free-flow travel times on all roads, which also corresponds to the fastest path. The route is computed by means of the Dynamic Programming Principle (DPP), where each road is assigned a weight equal to its estimated travel time. Once the route has been determined, the vehicle starts traveling towards its destination according to the second-order microscopic follow-the-leader model presented in Section \ref{sec:dynamics}.

During the simulation, whenever a compliant vehicle passes an RSU, it receives from the infrastructure the latest travel time estimates for all monitored roads. The OPU updates the corresponding road weights and recomputes the fastest path from the vehicle's current position to the assigned destination by means of the same DPP algorithm using the updated travel times. The newly computed path replaces the previous one from the current position onward, allowing the vehicle to adapt its route to the current traffic conditions. In contrast, non-compliant vehicles ignore the information provided by the infrastructure and continue following their initially computed route.

\subsubsection{Fuel consumption}
This traffic model also enables the estimation of the fuel consumption associated with the simulated traffic dynamics. Since the velocity and acceleration of every vehicle are available at each time step, the instantaneous fuel consumption rate can be computed according to the model proposed in \cite{borken2018comparing}. Assuming all vehicles to be gasoline average cars, the instantaneous fuel consumption rate $F_c(t)$ of vehicle $c$ is expressed as a quadratic function of the Vehicle Specific Power (VSP),
\begin{equation}\label{eq:fuel_cons}
F_c(t)=
\max\left\{
0,\,
m\left(
0.2102\,\mathrm{VSP}^2
+221\,\mathrm{VSP}
+596
\right)
\right\},
\end{equation}
where $m$ is the vehicle mass, set to $1.59$ ton. 
Following \citep{zhou2017development}, the VSP is estiamate as
\begin{equation}\label{eq:VSP}
\mathrm{VSP}
=
V_c\left(1.1A_c+9.81s+0.132\right)
+3.02\times10^{-4}V_c^3,
\end{equation}
where $V_c$ and $A_c$ denote the instantaneous velocity and acceleration of vehicle $c$, respectively, while $s$ is the road slope. In the considered simulations, flat roads are assumed, i.e., $s=0$. The coefficients are calibrated for $V_c$ expressed in m/s, $A_c$ in m/s$^2$, and $s$ in \%.
The model provides the instantaneous fuel consumption rate in g/h. This quantity is converted into L/s by dividing by $3600$ and by the gasoline specific weight, assumed to be $740$ g/L. The total fuel consumption is then obtained by integrating the instantaneous consumption rate of each vehicle over time and summing the contributions of all vehicles,
\begin{equation}\label{eq:total_fuel_cons}
J_F=\sum_{c=1}^{N_c}\int_{0}^{T}F_c(t)\,dt,
\end{equation}
where $T$ denotes the simulation horizon.

\subsection{RSUs optimization framework}
 Since the aim is to improve road network performance while minimizing the number of deployed RSUs, the methodology adopts a nested optimization scheme. In particular, an outer module explores the optimal number $k$ of RSUs, while an inner module, executed for each $k$, identifies the best placement configuration. The method is explained in the following with the detailed pseudocodes for all the algorithms provided in the \hyperref[sec:appendix]{Appendix}.
 \subsubsection{Placement optimization for fixed $k$}
Within the nested optimization framework, the RSUs placement optimization step acts as the inner loop, refining the configuration for each candidate number $k$ of RSUs. The search is driven by a multi-criteria method combined with a metaheuristic search.
The core logic is based on identifying high-impact road segments and iteratively refining their selection through simulation-based feedback. 
For a fixed number $k$ of RSUs, , the placement is computed through the \texttt{FindBestConfig}($k$) routine, formalized in \textbf{\hyperref[alg:rsu_optimizer]{Algorithm \ref*{alg:rsu_optimizer}}}.
This function works as the central block of the inner module that coordinates all the steps of the search. Its structure follows the standard Iterated Local Search (ILS) scheme described by \citet{lourencco2003iterated,lourencco2018iterated}, which is built around four basic elements: generating an initial solution, applying a local search, introducing a perturbation to escape local optima, and using an acceptance rule to decide how the search moves forward.
Specifically, the algorithm initializes the positions via \texttt{GetInitPos} function, evaluates the configuration $\mathbf{x}$ quality by computing $TTT_\mathbf{x}$ through traffic simulation (detailed in Section \ref{sec:traffic_simulation}), dynamically updates the link scores via \texttt{RankLinks}, and calls the \texttt{Perturbation} and \texttt{LocalSearch} routines whenever a stall condition (i.e., detection of a previously visited configuration $\mathbf{x} \in \mathcal{V}$) is met. The steps of the inner loop are detailed below, including initialization, road‑segment ranking, and the perturbation and local‑search steps. The procedure runs until the counter of consecutive non-improving iterations $c_{\text{fail}}$ reaches the maximum threshold $\tau_{\max}$. It finally returns, for the specific $k$, the best configuration $\mathbf{x}_k^*$ and the corresponding  $\Delta_{TTT}(\mathbf{x_k^*})$ to the high-level search strategy (either bisection or stepwise decrement). These strategies, described in detail in Section \ref{sec:k_opt}, differ in their search mechanism to determine the best $k$.

\paragraph{Initialization} \label{par:init}
The initial configuration $\mathbf{x}$ generated by $\texttt{GetInitPos}(G, M_{od}, k)$ function, as illustrated in \textbf{\hyperref[alg:rsu_optimizer]{Algorithm \ref*{alg:rsu_optimizer}}}, provides a deterministic starting point for the optimization process based on static shortest-path traffic flows. For every origin-destination vehicle request in $M_{od}$, the deterministic shortest path is computed using Dijkstra's algorithm,  taking the physical length of each link $l \in L$ as its weight. A counter is used to record the total number of shortest paths traversing each link $l$. The links are then sorted in descending order according to this path frequency counter. Finally, the initial RSU configuration $\mathbf{x}$ is constructed by assigning $x_l = 1$ to the top $k$ most-traversed road segments, and $x_l = 0$ to all remaining links.

\paragraph{Road segments ranking} \label{par:road_seg_rank}
The ranking of the network links  $L$, embedded in the $\texttt{RankLinks}(G, M_{od}, N_{\text{cars}}, P, \alpha)$ function in \textbf{\hyperref[alg:rsu_optimizer]{Algorithm \ref*{alg:rsu_optimizer}}}, is a fundamental step for generating candidate RSU configurations.  Each road segment $l$ is evaluated based on two distinct metrics: (i) \textit{path road coverage} $C_{path}(l)$ that counts the number of unique vehicle routes that traverse a specific road segment. By considering unique routes rather than raw vehicle counts, it identifies the importance of a link in terms of the number of different routes that traverse it. (ii) \textit{Link OD coverage} $C_{OD}(l)$ that quantifies for each origin-destination pair the number of distinct downstream path suffixes (i.e., the sequences of nodes leading from link $e$ to the destination $d$) that belong to the set of $k$-shortest paths. This metric identifies segments that act as critical decision points. A high $C_{OD}$ value indicates that a road segment might act as a diversion point to multiple routing alternatives. 

To prevent any bias introduced by scale differences, both metrics are normalized to the range $[0, 1]$ using min-max scaling. The final ranking score for each segment is computed as a weighted linear combination:
\begin{equation}
    \phi(l) = \alpha \cdot C_{OD}(l) + (1 - \alpha) \cdot C_{path}(l)
\end{equation}
In this study, we set $\alpha = 0.75$ to prioritize the topological strategic value $C_{OD}$ over the path-based coverage $C_{path}$ of each link $l$.
The resulting ranked list $L^r$ is then processed by the function $\texttt{SelectTopK}(L^r, k)$ to identify the $k$ best positions, which are used to generate the candidate configuration for RSUs placement.
\paragraph{Perturbation and Local search}
The perturbation step is designed to prevent the optimization process from getting trapped in local optima. This step is executed whenever the ranking step produces a configuration $\mathbf{x}$ that has already been evaluated, as recorded in the set of visited solutions.  The procedure detailed in \textbf{\hyperref[alg:perturbation]{Algorithm \ref*{alg:perturbation}}} operates by identifying a specific subset of $m$ RSUs, denoted as $R_{\mathbf{x}}^{mobile}$, which are relocated to different road segments while strictly preserving the total number of RSUs $k$. These $m$ units are removed from their current positions and reassigned to $m$ random positions within the set of empty segments $L_{\mathbf{x}}^{empty}$. The process iterates until a unique configuration $\mathbf{x}^{pert} \notin \mathcal{V}$ is found.

Once a unique solution is identified through perturbation, a local search is performed, as detailed in \textbf{\hyperref[alg:localSearch]{Algorithm \ref*{alg:localSearch}}}. The algorithm examines the road network to identify the set of adjacent segments $Adj(u)$ for each mobile RSU $u \in R_{\mathbf{x}}^{mobile}$, defined as links sharing a common start or end node. A set of candidate configurations $\mathcal{N}$ is generated by systematically moving the mobile RSUs to these adjacent links while maintaining the positions of the RSUs in $R_{\mathbf{x}}^{fixed}$. Every configuration in this neighborhood is evaluated via the traffic simulation. The configuration $\mathbf{x}$ that yields the largest $TTT$ reduction compared to the baseline, i.e., the maximum $\Delta_{TTT}$ in the neighborhood, becomes the new reference configuration for the next iteration, even if it does not improve upon the previously evaluated solution, in order to avoid local optima.
\subsubsection{$k$-optimization} \label{sec:k_opt}
To determine the optimal number of devices, the outer loop relies on two alternative strategies that manage the RSUs number $k$ differently: the \textit{Stepwise Decrement Approach} (\texttt{SD}) and the \textit{Bisection Search} (\texttt{BS}). 

The \texttt{SD} approach presented in \textbf{\hyperref[alg:stepwise]{Algorithm \ref*{alg:stepwise}}} explores the solution space by starting from a maximum number of RSUs and iteratively reducing the number of units. If the local search fails to improve the score after a threshold of $\tau_{\max}$ iterations, the number of RSUs $k$ is decremented of $\kappa$ units.

The \texttt{BS} approach described in \textbf{\hyperref[alg:bisection]{Algorithm \ref*{alg:bisection}}} optimizes the number of RSUs $k$ by treating it as a partition problem. It evaluates the boundary conditions ($k_{min}$, $k_{max}$) and the midpoint, then narrows the interval toward the most promising value of $k^*$ and the best configuration $\mathbf{x}^*$ given that value of $k$. Both the approaches terminate when a maximum number of iterations $I_{\max}$ is reached, returning the best configuration $\mathbf{x}^*$.

\section{Results} \label{sec:results}

\subsection{Geometries }
\begin{figure}[H] 
	\centering  \includegraphics[width=0.75\linewidth]{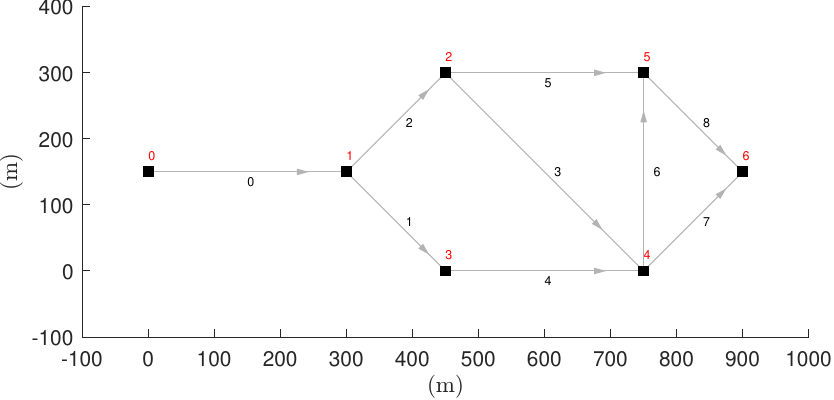}
    \caption{Diamond network with road and junction numbering.}
    \label{fig:Diamond_network}
\end{figure}

\begin{figure}[H] 
	\centering      \includegraphics[width=0.6\linewidth]{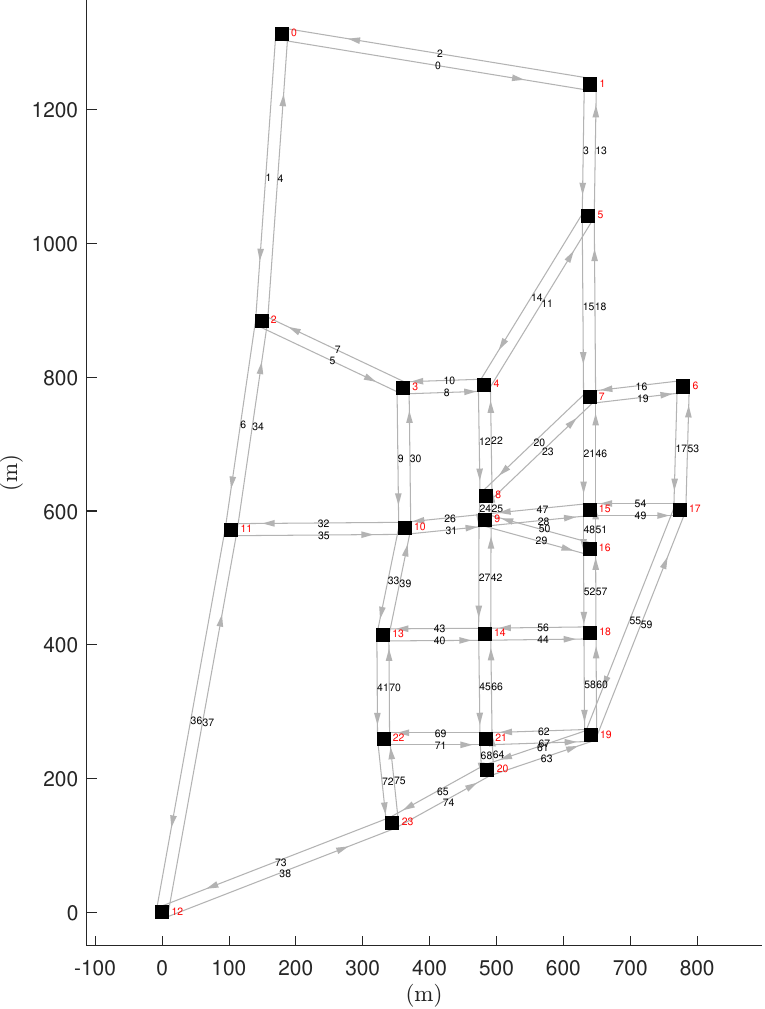}
    \caption{Sioux Falls network with road and junction numbering.}
    \label{fig:SiouxFalls_network}
\end{figure}

To evaluate the effectiveness of the proposed methods, we first considered a small-scale baseline network referred to as the Diamond network because of its geometric shape. The compact size of this network, which consists of 7 nodes and 9 links (as depicted in Fig.\ \ref{fig:Diamond_network}) was intentionally selected to enable the execution of an exhaustive search (\texttt{ES}), which is defined in \textbf{\hyperref[alg:exhaustive]{Algorithm \ref*{alg:exhaustive}}}.
As an additional network for our tests, we considered Sioux Falls, a real-world instance\footnote{Transportation networks for research github repository, available at \url{https://github.com/bstabler/TransportationNetworks}} (see \citealp{stabler2018transportation}) consisting of 24 nodes and 76 links (as illustrated in Fig.\  \ref{fig:SiouxFalls_network}). Networks characteristics are summarized in Table \ref{tab:net_params}.
\begin{table}[H]
    \centering
    \caption{Networks sizes}
    \label{tab:net_params}
    \resizebox{0.3\textwidth}{!}{
    \begin{tabular}{lcc}
        \toprule
        Network & Nodes & Links   \\
        \midrule
         Diamond & 7 & 9 \\
         Sioux Falls & 24 & 76 \\
        \bottomrule
    \end{tabular}
    }
\end{table}
\vfill
\newpage


\subsection{Computational setup}
The proposed methodology has been implemented using a serial code, with the optimization module written in Python and the dynamic simulation developed in C. All experiments were executed on a server equipped with an Intel Xeon Silver 4214R (2.40 GHz) and 128 GB RAM.

The microscopic traffic model is discretized in time by means of a semi-implicit Euler scheme with a fixed time step $\Delta t$. At each time step, the acceleration of each vehicle is first computed from \eqref{eq:acc}. The vehicle velocity is then updated, followed by the position using the newly computed velocity. After each update, it is verified whether the vehicle has reached the end of its current road. If so, the vehicle is transferred to the next road along its route or removed from the simulation when its destination is reached.
The model parameters are kept fixed throughout all tests and are reported in Table \ref{tab:traffic_params}.

\begin{table}[htbp]
    \centering
    \caption{Parameters used in numerical simulations.}
    \label{tab:traffic_params}
    \resizebox{0.5\textwidth}{!}{
    \begin{tabular}{lccc}
        \toprule
        Parameter & Symbol & Value & Unit \\
        \midrule
        Time step & $\Delta t$ & 0.6 & s\\
        Maximum admissible velocity & $V_{\max}$ & 50.0 & km/h \\
        Minimum safe distance & $\delta_{\mathrm{close}}$ & 10.0 & m \\
        Free-flow distance threshold & $\delta_{\mathrm{far}}$ & 40.0 & m \\
        Acceleration time & $\tau_{\mathrm{acc}}$ & 6 & s \\
        Deceleration time & $\tau_{\mathrm{dec}}$ & 0.6 & s \\
        \bottomrule
    \end{tabular}
    }
\end{table}
With respect to the parameter setting in the optimization module, \texttt{the SD} and \texttt{BS} strategies share the same parameter values for $I_{max}$ (the maximum number of iterations) and $\tau_{max}$ (the maximum number of consecutive non-‐improving iterations). For the \texttt{SD} method, we additionally defined the integer decrement step $\kappa$. These parameters were determined through a numerical iterative procedure aimed at determining proper values for each network, and are summarized in Table \ref{tab:algo_params}.

\begin{table}[htbp]
    \centering
    \caption{Algorithms parameters}
    \label{tab:algo_params}
    \resizebox{0.3\textwidth}{!}{
    \begin{tabular}{lccc}
        \toprule
        Network & $I_{max}$ & $\tau_{max}$ & $\kappa$  \\
        \midrule
         Diamond & 50 & 5 & 1  \\
         Sioux Falls & 150 & 10 & 5 \\
        \bottomrule
    \end{tabular}
    }
\end{table}

\subsection{Why RSUs placement is a complex task}\label{sec:trade-off}
To begin with, we propose a numerical test showing why nontrivial optimization algorithms are needed for RSUs placement. 
The complexity mainly comes from the fact that the information collected by RSUs depends on the routing decisions, which, in turn, are (at least partially) generated by the RSUs themselves. 
Therefore, optimizing RSUs placement requires accounting for its feedback on re-routing.
A naive strategy would be to deploy an RSU on a road with high traffic flow, since a large number of passing vehicles enables the RSU to continuously collect a large amount of up-to-date information and, therefore, accurately estimate the travel time of the monitored road. 
However, this choice may be counterproductive. 
Once the RSUs detect high car flux or congestion, they start recommending avoiding the monitored road, moving vehicles to roads that are possibly not monitored. 
As a result, the amount of information available to update travel times is reduced, path suggestions can become inaccurate, and vehicles could be routed towards longer or more congested paths.

This highlights the main challenge in RSUs placement, which requires balancing two conflicting objectives: on the one hand, RSUs should monitor roads with high traffic flow to collect as much information as possible and accurately estimate travel times; on the other hand, they should spread vehicles across the whole road network to better exploit road usage and reduce congestion. 
Therefore, an effective RSUs placement must strike a balance between information collection and traffic redistribution.
Interestingly, a similar situation is found in \cite{cristiani2026microscopic}, where vehicles paradoxically need to encounter as many other vehicles as possible in order to reduce the chance of ending up in a traffic jam. 

To illustrate this trade-off, we consider a scenario on the Sioux Falls network (Fig.\ \ref{fig:SiouxFalls_network}) with 150 vehicles: 100 vehicles start from junction 19, 30 from junction 13, and 20 from junction 12, all having junction 0 as their destination. 
Initially, all vehicles follow the shortest path to the destination. 
We assume full compliance, therefore whenever a vehicle encounters an RSU, it follows the routing recommendation computed using the information currently available at that RSU.

\begin{figure}[H]
    \centering
    \subfloat[reference configuration (with two RSUs) (\href{www.emilianocristiani.it/attach/paper_I2X/siouxfalls_2RSU_in26-50.mp4}{video})]{
    \label{subfig:ref_conf}
    \includegraphics[width=.28\linewidth]{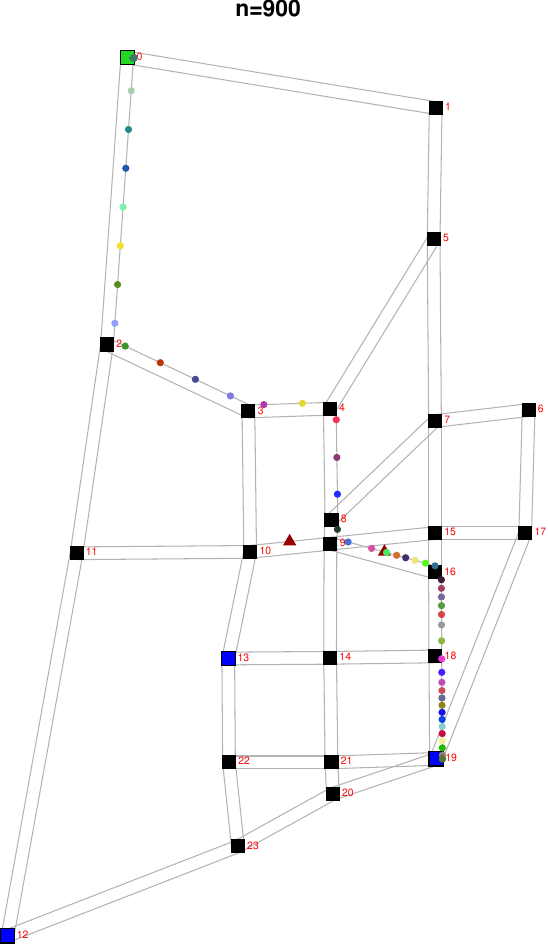}}\\
    \subfloat[third RSU added on road 7 
    (\href{www.emilianocristiani.it/attach/paper_I2X/siouxfalls_3RSU_in7-26-50.mp4}{video})]{
    \label{subfig:worst_conf}
    \includegraphics[width=.28\linewidth]{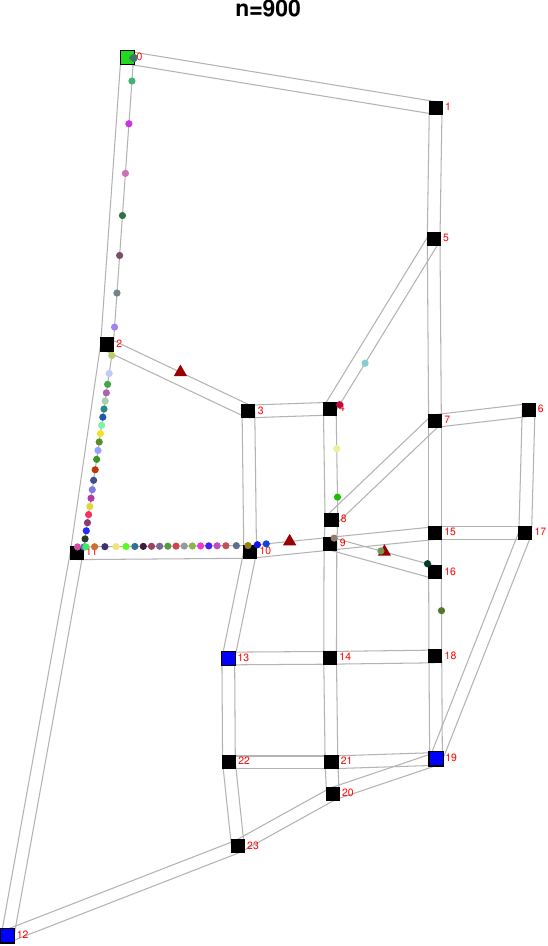}}\quad
    \subfloat[third RSU added on road 10 
    (\href{www.emilianocristiani.it/attach/paper_I2X/siouxfalls_3RSU_in10-26-50.mp4}{video})]{
    \label{subfig:best_conf}
    \includegraphics[width=.28\linewidth]{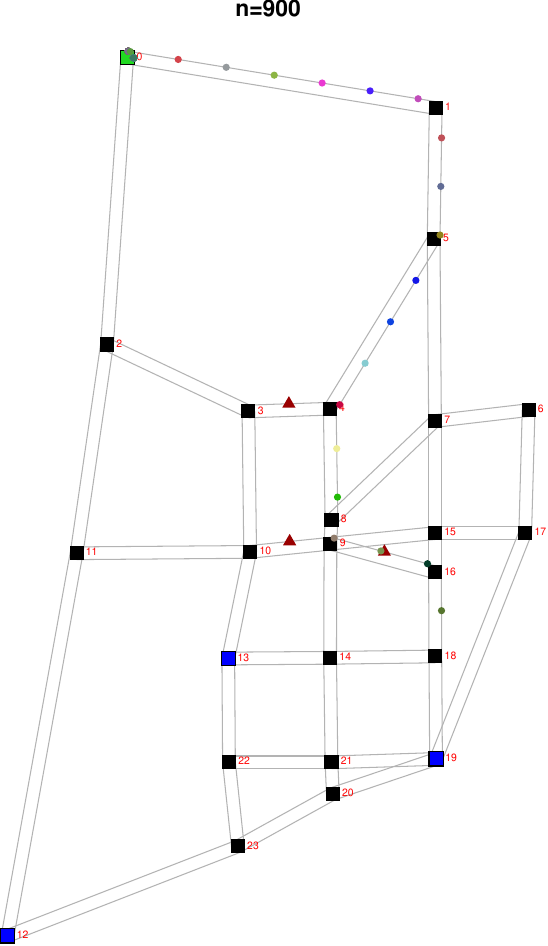}}\par
    \caption{Sioux Falls network: three screenshots illustrating different RSUs placement configurations, all taken at the same simulation time. The comparison highlights the different traffic conditions resulting from the three RSUs placements. Blue junctions indicate vehicle origins, the green junction indicates the common destination, and red triangles represent the RSUs.}
    \label{fig:trade-off}
\end{figure}

As a reference configuration, we place two RSUs on roads 26 and 50 (Fig.\ \ref{subfig:ref_conf}). Under this configuration, only vehicles departing from junction 19 receive routing recommendations, and the resulting $TTT$ is $1058.93$ minutes. Excluding road 4, which naturally carries the highest traffic because it is one of the final roads leading to the common destination, road 7 has the highest traffic flow in the network. 
We then place a third RSU on road 7, as shown in Fig.\ \ref{subfig:worst_conf}. Although this appears to be a natural choice, the $TTT$ increases to $1102.28$ minutes. 
Once congestion develops on road 7, the RSU redirects vehicles towards an alternative, partially nonmonitored, longer route. 
This initially relieves congestion on road 7. 
However, after some time, even when traffic conditions on road 7 improve, vehicles continue to be rerouted to the longer alternative route because the RSU does not receive enough information to immediately update the travel time of road 7. 

As an alternative strategy, we place the third RSU on road 10 (Fig.\ \ref{subfig:best_conf}), which carries a moderate traffic flow in the reference configuration. As in the previous configuration, after some time the RSU on road 10 also stops receiving updated information. Nevertheless, its location allows it to redistribute traffic across the available routes before severe congestion develops on road 7, preventing the formation of a bottleneck and reducing the overall $TTT$ to $875.81$ minutes.

We can conclude that simply increasing the number of RSUs or deploying them on the roads with the highest traffic flow does not necessarily improve network performance. The effectiveness of an RSU depends not only on the amount of information it can collect but also on how its routing recommendations modify future traffic and, consequently, the information available for future routing decisions. 

\subsection{Diamond network}
\begin{figure}[H] 
    \centering 
    \subfloat[without RSUs (\href{www.emilianocristiani.it/attach/paper_I2X/diamond_0RSU.mp4}{video})]{
    \label{subfig:0RSU}
    \includegraphics[width=.45\linewidth]{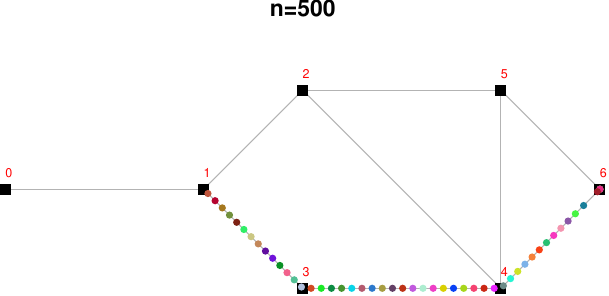}}\quad 
    \subfloat[optimal RSU placement (\href{www.emilianocristiani.it/attach/paper_I2X/diamond_ES_opt_conf_2RSU.mp4}{video})]{ \label{subfig:opt_conf} \includegraphics[width=.45\linewidth]{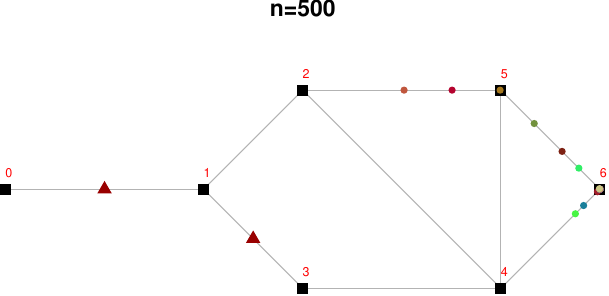}} 
    \caption{Diamond network: comparison between the case without RSUs and the optimal RSUs placement identified by the \texttt{ES}. The screenshots are taken at the same simulation time, allowing a direct comparison of the resulting traffic conditions.} 
    \label{fig:opt_conf_ES} 
\end{figure}

By computing the global optimum through the \texttt{ES}, we could verify whether the \texttt{BS} and the \texttt{SD} strategies were able to converge to the optimal solution. Due to the small scale of the network, we set $\kappa = 1$ for the stepwise decrement approach, thereby decreasing $k$ in steps of one. In these preliminary numerical tests, two key parameters were systematically varied to assess their impact on algorithm performance. First, we analyzed the compliance rate, denoted by $\gamma$ and defined as the proportion of users who actively adhere to the provided rerouting strategies; this parameter was tested at four levels, specifically $\gamma \in \{0.25, 0.50, 0.75, 1.00 \}$. Second, we varied the total number of vehicles within the network, represented by $N_c$. This variable represents the aggregate user demand within the reference time frame and was varied across five different traffic levels: $N_c \in \{ 50, 75, 100, 125, 150 \}$. By changing $N_c$ in this way, we were able to evaluate how different levels of network congestion influence both the computational times and the quality of the solutions generated by the proposed algorithms. In the considered scenario, all vehicles have junction 6 as their destination. Most vehicles depart from junction 0, while two vehicles depart from junctions 1 and 2, respectively. These two vehicles travel at constant speeds of 9 km/h and 12 km/h, respectively, thereby causing a general slowdown of the traffic flow throughout the network.
\begin{figure}[H]
    \centering
    \includegraphics[width=0.65\linewidth]{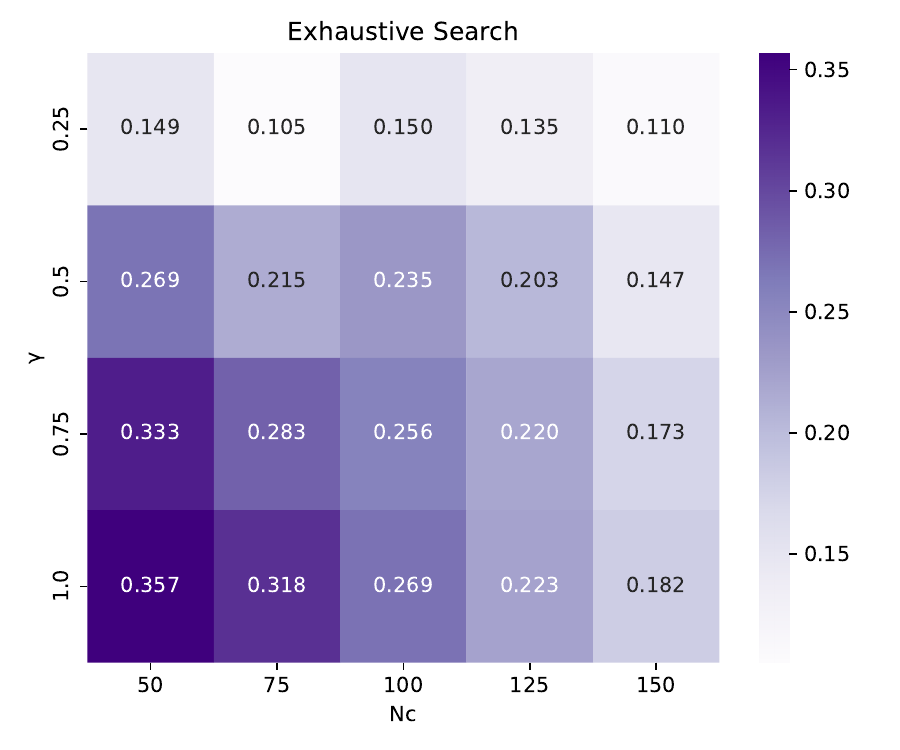}
    \caption{Diamond network: $\Gamma_{TTT}$ resulting from exhaustive search as compliance $\gamma$ varies and for different demand levels $N_{c}$}
    \label{fig:heatmaps}
\end{figure}
The heatmaps depicted in Fig.\ \ref{fig:heatmaps} show the $\Gamma_{TTT}$, which represents the percentage gap between the $TTT_{\max}$ and $TTT_{\min}$  found by the \texttt{ES}. This visually highlights how changes in network demand $N_c$ and user compliance $\gamma$ actually impact the performance. In particular, for any given compliance level $\gamma$, an increase in network demand $N_c$ naturally induces a generalized increase of both $TTT_{\min}$ and $TTT_{\max}$ due to physical capacity constraints. However, the percentage gap $\Delta_{TTT}$ drops sharply as the system gets heavily congested. For instance, at $\gamma = 1.0$, $\Delta_{TTT}$ falls from over $0.36$ down to $0.18$. This behaviour suggests that, once the network is congested, the system becomes less responsive to traffic control strategies, as a saturated network leaves little room for any strategy to deliver significant improvements. On the other hand, increasing the compliance rate significantly amplifies the impact of traffic control policies, expanding the gap between the best and worst-case network states. At a low compliance level of $\gamma = 0.25$, the selfish driver behavior severely limits the algorithms effectiveness, keeping $\Delta_{TTT}$ confined to the 0.10–0.15 range. This confirms that strategic traffic control is most powerful when users fully adhere to the route guidance directives. Fig.\ \ref{fig:opt_conf_ES} compares the traffic conditions obtained without RSUs and with the optimal RSUs placement identified by the \texttt{ES}, for $N_c = 50$ and $\gamma = 1.00$. The optimal configuration consists of two RSUs deployed on roads 0 and 1 (Fig.\ \ref{subfig:opt_conf}). 

An examination of the $TTT_{\min}$ values reported in Table \ref{tab:ttt_results_diamond} reveals consistency among the three approaches in terms of solution quality. Regardless of the demand level or user compliance rate, all three approaches systematically converge to the same $TTT_{\min}$. 
Interestingly, we also found it valuable to compute the worst-case network states $TTT_{\max}$ through the \texttt{ES}. This analysis revealed that the maximum delay does not simply coincide with the baseline scenario i.e., $k = 0$. Instead, $k_{\max}$ typically falls around higher control thresholds, such as $k = 6$ or $k = 7$. Thus, network delay can remain high even as $k$ increases, showing that the reduction in $TTT$ is not linear with the increase in RSUs. This makes the problem far from trivial.
\begin{table}[H]
    \centering
    \caption{Diamond network: $(TTT, k)$, $\Gamma_{TTT}$ and $\Delta_{TTT}$ results as compliance $\gamma$ varies and for different demand levels $N_c$. $TTT$ is expressed in minutes.}
    \label{tab:ttt_results_diamond}
    \small
    \setlength{\tabcolsep}{3pt} 
    \resizebox{\textwidth}{!}{
    \begin{tabular}{cc ccc ccc ccc}
        \toprule
        & & \multicolumn{3}{c}{\texttt{ES}} & \multicolumn{3}{c}{\texttt{BS}} & \multicolumn{3}{c}{\texttt{SD}} \\
        \cmidrule(lr){3-5} \cmidrule(lr){6-8} \cmidrule(lr){9-11}
        $\gamma$ & $N_c$ & $(TTT_{\min}, k_{\min})$ & $(TTT_{\max}, k_{\max})$ & $\Gamma_{TTT} (\%)$ & $(TTT_{\min}, k_{\min})$ & $TTT_0$ & $\Delta_{TTT} (\%)$  & $(TTT_{\min}, k_{\min})$ & $TTT_0$ & $\Delta_{TTT}$ (\%)\\
        \midrule
        0.25 & 50  & (270.16, 2) & (318.25, 6) & 15.11 & (270.16, 4) & 317.52 & 14.92 & (270.16, 5) & 317.52 & 14.92 \\
             & 75  & (482.76, 2) & (540.31, 6) & 10.65 & (482.76, 4) & 539.48 & 10.51 & (482.76, 3) & 539.48 & 10.51 \\
             & 100 & (682.00, 2) & (812.37, 7) & 16.05 & (682.00, 4) & 802.35 & 15.00 & (682.00, 4) & 802.35 & 15.00 \\
             & 125 & (956.95, 2) & (1114.13, 7)& 14.11 & (956.95, 4) & 1106.03& 13.48 & (956.95, 5) & 1106.03& 13.48 \\
             & 150 & (1291.22, 2)& (1473.08, 7)& 12.35 & (1291.22, 4)& 1450.52& 10.98 & (1291.22, 3)& 1450.52& 10.98 \\
        \midrule
        0.50 & 50  & (232.05, 2) & (319.96, 6) & 27.48 & (232.05, 4) & 317.52 & 26.92 & (232.05, 3) & 317.52 & 26.92 \\
             & 75  & (423.36, 2) & (548.35, 7) & 22.79 & (423.36, 4) & 539.48 & 21.52 & (423.36, 4) & 539.48 & 21.52 \\
             & 100 & (613.93, 2) & (841.89, 7) & 27.08 & (613.93, 4) & 802.35 & 23.48 & (613.93, 4) & 802.35 & 23.48 \\
             & 125 & (880.96, 2) & (1167.08, 7)& 24.52 & (880.96, 4) & 1106.03& 20.35 & (880.96, 4) & 1106.03& 20.35 \\
             & 150 & (1237.28, 2)& (1567.22, 7)& 21.05 & (1237.28, 4)& 1450.52& 14.70 & (1237.28, 4)& 1450.52& 14.70 \\
        \midrule
        0.75 & 50  & (211.82, 2) & (321.14, 6) & 34.04 & (211.82, 4) & 317.52 & 33.29 & (211.82, 3) & 317.52 & 33.29 \\
             & 75  & (386.91, 2) & (550.14, 7) & 29.67 & (386.91, 4) & 539.48 & 28.28 & (386.91, 4) & 539.48 & 28.28 \\
             & 100 & (596.95, 2) & (844.19, 6) & 29.29 & (596.95, 4) & 802.35 & 25.60 & (596.95, 4) & 802.35 & 25.60 \\
             & 125 & (862.49, 2) & (1176.41, 6)& 26.68 & (862.49, 4) & 1106.03& 22.02 & (862.49, 3) & 1106.03& 22.02 \\
             & 150 & (1199.00, 2)& (1565.62, 6)& 23.42 & (1199.00, 4)& 1450.52& 17.34 & (1199.00, 4)& 1450.52& 17.34 \\
        \midrule
        1.00 & 50  & (204.21, 2) & (322.68, 6) & 36.71 & (204.21, 4) & 317.52 & 35.69 & (204.21, 3) & 317.52 & 35.69 \\
             & 75  & (367.94, 2) & (545.24, 6) & 32.52 & (367.94, 4) & 539.48 & 31.79 & (367.94, 4) & 539.48 & 31.79 \\
             & 100 & (586.18, 2) & (815.69, 6) & 28.14 & (586.18, 4) & 802.35 & 26.94 & (586.18, 7) & 802.35 & 26.94 \\
             & 125 & (858.94, 2) & (1128.51, 6)& 23.89 & (858.94, 4) & 1106.03& 22.34 & (858.94, 4) & 1106.03& 22.34 \\
             & 150 & (1186.21, 2)& (1483.64, 6)& 20.05 & (1186.21, 4)& 1450.52& 18.22 & (1186.21, 4)& 1450.52& 18.22 \\
        \bottomrule
    \end{tabular}
    }
\end{table}
The execution times reported in Table \ref{tab:execution_times} show significant differences in computational effort, reflecting the trade-off between accuracy and efficiency. \texttt{ES} serves as the global‑optimality baseline, but its runtime grows rapidly with $N_c$, increasing from about 19 to nearly 78 minutes. This makes it unsuitable for large‑scale networks. The stepwise decrement strategy represents a middle-ground option, by cutting execution times by roughly $50\%$ to $60\%$ compared to the  \texttt{ES} (e.g., taking $34.8$ minutes instead of $77.5$ minutes at $\gamma = 1.0$ and $N_c = 150$). However, its execution time remains relatively high because it increments sequentially through the solution space. In this context, the \texttt{BS} demonstrates computational superiority by achieving an order-of-magnitude reduction in processing times. Across all tested scenarios, the bisection framework restricts the runtimes to well under a $9$ minutes threshold, operating between $5\times$ to $10\times$ faster than the stepwise strategy, running up to $13\times$ faster than the exhaustive baseline, and then showing scalability.
\begin{table}[htbp]
    \centering
    \caption{Diamond network: computational times (sec) as compliance $\gamma$ varies and for different demand levels $N_{c}$}
    \label{tab:execution_times}
    \small
    \resizebox{0.55\textwidth}{!}{
    \begin{tabular}{lcccc || ccccc}
        \toprule
        $\gamma$ & $N_{c}$ & \texttt{ES} & \texttt{BS}  & \texttt{SD} & $\gamma$ & $N_{c}$ & \texttt{ES} & \texttt{BS}  & \texttt{SD} \\
        \midrule
        0.25 &  50  & 19 & 2 &  9 & 0.75 &  50  & 19 & 2 &  8 \\
        0.25 &  75  & 30 & 2 & 14 & 0.75 &  75  & 31 & 4 & 14 \\
        0.25 & 100  & 44 & 4 & 21 & 0.75 & 100  & 44 & 5 & 18 \\
        0.25 & 125  & 60 & 5 & 26 & 0.75 & 125  & 60 & 5 & 29 \\
        0.25 & 150  & 78 & 6 & 36 & 0.75 & 150  & 78 & 9 & 34 \\
        \midrule
        0.50 &  50  & 20 & 2 &  9 & 1.00 &  50  & 19 & 2 &  9 \\
        0.50 &  75  & 30 & 2 & 13 & 1.00 &  75  & 29 & 3 & 11 \\
        0.50 & 100  & 45 & 5 & 20 & 1.00 & 100  & 43 & 5 & 19 \\
        0.50 & 125  & 60 & 5 & 27 & 1.00 & 125  & 60 & 6 & 26 \\
        0.50 & 150  & 75 & 8 & 33 & 1.00 & 150  & 78 & 9 & 35 \\
        \bottomrule
    \end{tabular}
    }
\end{table}
\subsection{Sioux Falls network}
After validating the proposed methodology on the Diamond network, we evaluate its performance on the Sioux Falls network, which represents a considerably more challenging benchmark due to its larger size and the presence of multiple alternative routes between origin-destination (OD) pairs.

In what follows, we present the results for the Sioux Falls network in a \emph{deterministic case} (Section \ref{sec:siouxfalls_det}), where we consider two demand scenarios, namely S1 and S2, and in a \emph{stochastic case} (Section \ref{sec:siouxfalls_stoc}). 
In both cases, we set the number of vehicles $N_c$ to 550, a value obtained from a calibration process aimed at determining a reasonable traffic volume. 
For this network, only bisection search (\texttt{BS}) and stepwise decrement (\texttt{SD}) methods were applied, as performing the exhaustive search (\texttt{ES}) was computationally intractable. 

\subsubsection{Sioux Falls: deterministic demand case}
\label{sec:siouxfalls_det}
In this section, we consider two deterministic traffic demand scenarios characterized by fixed OD pairs.
The selected demand patterns are inspired by the enriched Sioux Falls scenario proposed by \cite{chakirov2014enriched}, who extended the classical Sioux Falls benchmark by integrating real demographic information, land-use data, census data, building information, and a synthetic population model. Their framework provides a realistic characterization of home and work locations within the Sioux Falls network.

\begin{figure}[H] 
	\centering  
    \includegraphics[width=0.4\linewidth]{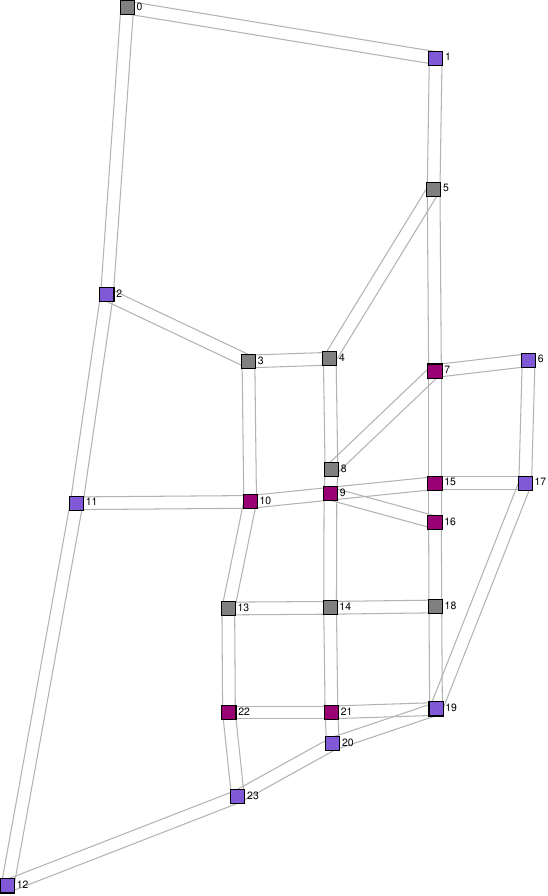}
    \caption{Sioux Falls network: classification of the junctions. Residential junctions are shown in purple, work junctions in magenta, and mixed junctions in gray.}
    \label{fig:junc_class}
\end{figure}

Based on this characterization, and as illustrated in Fig.\ \ref{fig:junc_class}, we classify the network junctions into three categories: residential, work, and mixed. Residential junctions correspond to areas with a high concentration of home locations, work junctions correspond to areas with a high concentration of workplaces, whereas mixed junctions represent intermediate areas where neither residential nor employment activities clearly prevail.

Using this classification, we define the two deterministic traffic demand scenarios. In scenario S1, vehicles depart from residential junctions and travel towards work junctions, reproducing the predominant home-to-work commuting pattern during the morning peak period. Conversely, scenario S2 considers the opposite traffic demand, with vehicles traveling from work junctions towards residential junctions. In particular, we consider the OD pairs reported in Table \ref{tab:deterministic_OD}. Vehicles sharing the same origin are released sequentially through a delayed departure mechanism. This prevents vehicles approaching the origin junction from other incoming roads from being forced to wait until all vehicles starting at that junction have entered the network, thus avoiding unrealistic congestion effects at the beginning of the simulation. 

\begin{table}[ht]
\centering
\caption{Sioux Falls network: Deterministic OD demand scenarios. Both scenarios consist of 550 vehicles, with 50 vehicles assigned to each OD pair.}
\label{tab:deterministic_OD}
\resizebox{0.65\textwidth}{!}{
\begin{tabular}{|c|p{9cm}|}
\hline
\textbf{Scenario} & \textbf{OD pairs} \\
\hline
S1 &
$(12,15)$, $(1,15)$, $(1,10)$, $(6,10)$, $(6,21)$,
$(17,22)$, $(2,22)$, $(11,22)$, $(20,7)$,
$(19,9)$, $(23,16)$
\\
\hline
S2 &
$(15,12)$, $(15,1)$, $(10,1)$, $(10,6)$, $(21,6)$,
$(22,17)$, $(22,2)$, $(22,11)$, $(7,20)$,
$(9,19)$, $(16,23)$
\\
\hline
\end{tabular}
}
\end{table}

\begin{figure}[H]
    \centering
    \subfloat[scenario S1 without RSUs\\ \centering (\href{www.emilianocristiani.it/attach/paper_I2X/siouxfalls_S1_0RSU.mp4}{video})]{
    \label{subfig:S1_0RSU}
    \includegraphics[width=.33\linewidth]{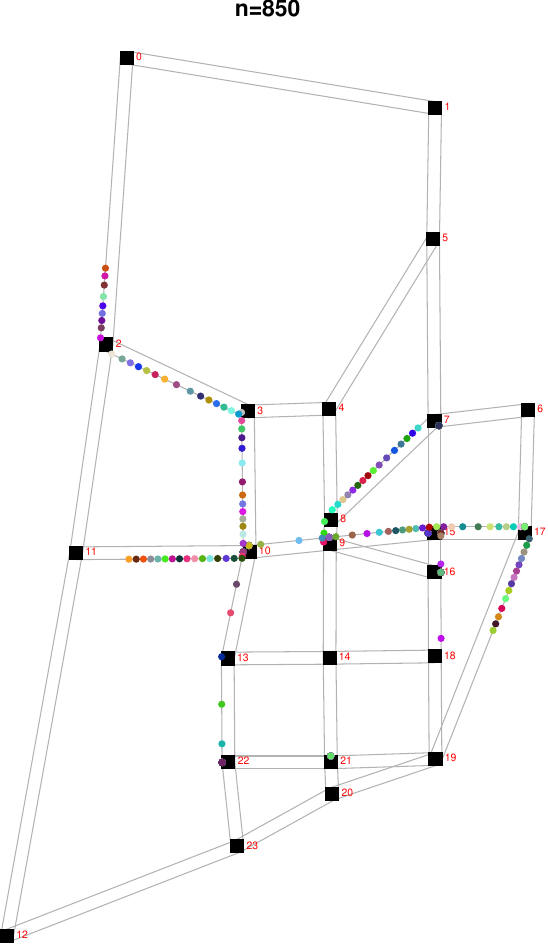}}\quad
    \subfloat[scenario S1 with the optimal RSUs placement (\href{www.emilianocristiani.it/attach/paper_I2X/siouxfalls_S1_opt_conf_21RSU.mp4}{video})]{
    \label{subfig:S1_opt_conf}
    \includegraphics[width=.33\linewidth]{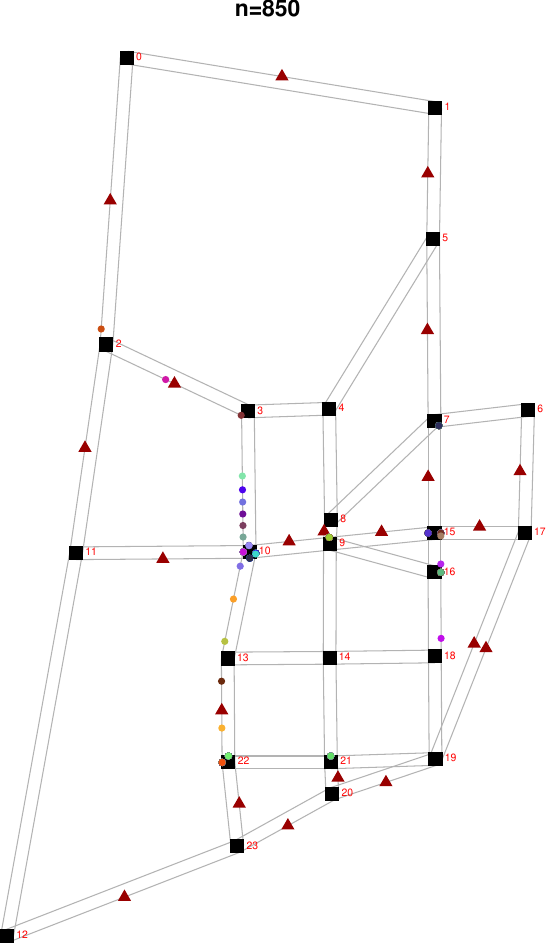}}\\
    \subfloat[scenario S2 without RSUs \\
    \centering (\href{www.emilianocristiani.it/attach/paper_I2X/siouxfalls_S2_0RSU.mp4}{video})]{
    \label{subfig:S2_0RSU}
    \includegraphics[width=.33\linewidth]{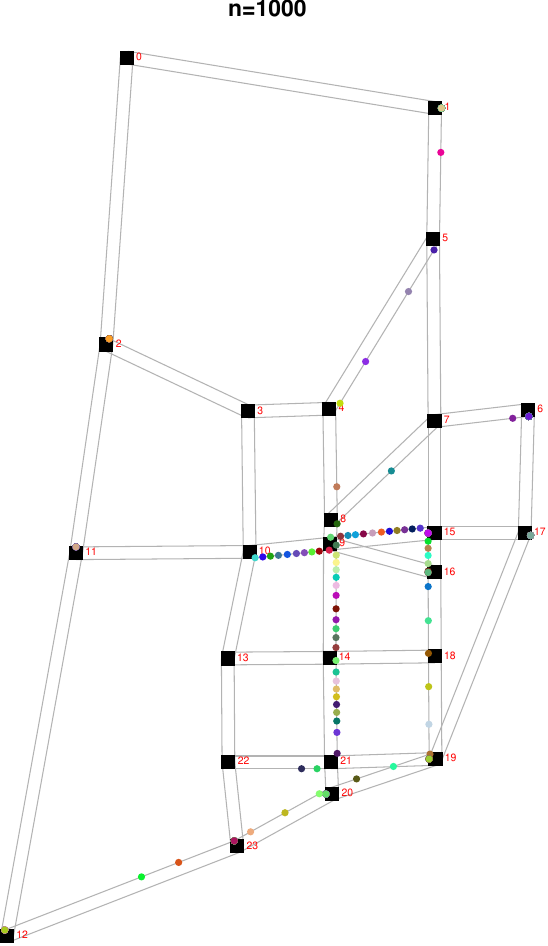}}\quad
    \subfloat[scenario S2 with the optimal RSUs placement (\href{www.emilianocristiani.it/attach/paper_I2X/siouxfalls_S2_opt_conf_66RSU.mp4}{video})]{
    \label{subfig:S2_opt_conf}
    \includegraphics[width=.33\linewidth]{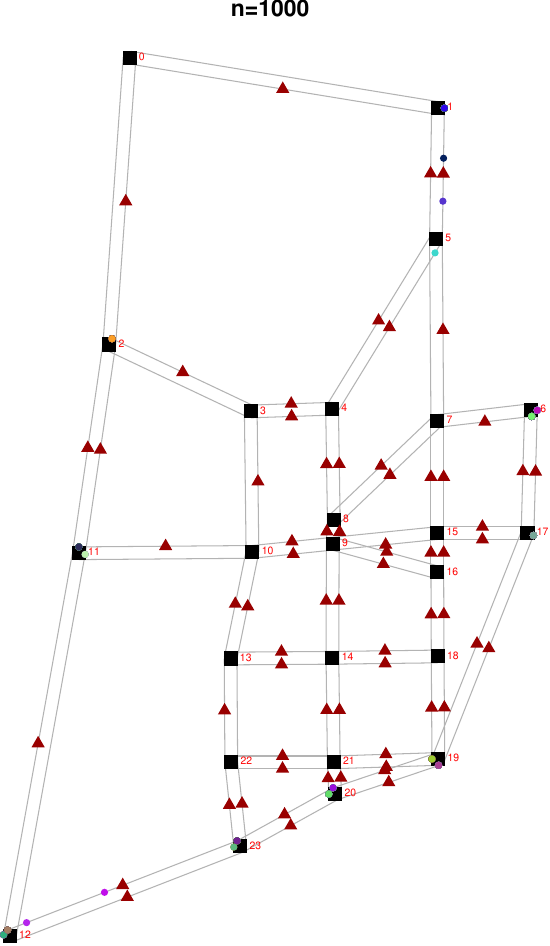}}\quad
    \caption{Sioux Falls network: comparison between the case without RSUs and the optimal RSUs placement for deterministic scenarios S1 and S2. For each scenario, the screenshots correspond to the same simulation time, enabling a comparison of the resulting traffic conditions.}
    \label{fig:opt_conf_determ}
\end{figure}

Tables \ref{tab:ttt_results_siouxfalls_s1} and \ref{tab:ttt_results_siouxfalls_s2} show the $TTT_{\min}$ results for the S1 and S2 scenarios using both the \texttt{BS} and \texttt{SD} methods, evaluated across compliance levels $\gamma \in \{0.25, 0.50, 0.75, 1.00\}$ and a fixed demand level of $N_c = 550$. The computational times associated with both methods, for scenarios S1 and S2 are presented in Tables \ref{tab:execution_times_s1} and \ref{tab:execution_times_s2}, respectively.
\begin{table}[H]
    \centering
    \caption{Sioux Falls network, deterministic case (Scenario S1): $(TTT, k)$ and $\Delta_{TTT}$ results as compliance $\gamma$ varies for demand level $N_c = 550$. $TTT$ is expressed in minutes.}
    \label{tab:ttt_results_siouxfalls_s1}
    \resizebox{0.8\textwidth}{!}{
    \begin{tabular}{cc ccc ccc}
        \toprule
        & & \multicolumn{3}{c}{\texttt{BS}} & \multicolumn{3}{c}{\texttt{SD}} \\
        \cmidrule(lr){3-5} \cmidrule(lr){6-8}
        $\gamma$ & $N_c$ & $(TTT_{\min}, k_{\min})$ & $TTT_0$ & $\Delta_{TTT}$ (\%) & $(TTT_{\min}, k_{\min})$ & $TTT_0$ & $\Delta_{TTT}$ (\%) \\
        \midrule
        0.25 & 550 & (2944.55, 38) & 3352.46 & 12.17 & (2928.34, 26) & 3352.46 & 12.65 \\
        0.50 & 550 & (2789.18, 47) & 3352.46 & 16.80 & (2789.10, 41) & 3352.46 & 16.80 \\
        0.75 & 550 & (2740.31, 38) & 3352.46 & 18.26 & (2717.54, 31) & 3352.46 & 18.94 \\
        1.00 & 550 & (2685.32, 42) & 3352.46 & 19.90 & (2674.62, 21) & 3352.46 & 20.22 \\
        \bottomrule
    \end{tabular}
    }
\end{table}
\begin{table}[H]
    \centering
    \caption{Sioux Falls network, deterministic case (Scenario S1): computational times (sec) as compliance $\gamma$ varies for demand level $N_{c} = 550$}
    \label{tab:execution_times_s1}
    \resizebox{0.25\textwidth}{!}{
    \begin{tabular}{cccc}
        \toprule
        $\gamma$ & $N_{c}$ & \texttt{BS}  & \texttt{SD} \\
        \midrule
        0.25 & 550 & 700 & 1102 \\
        0.50 & 550 & 504 & 1088 \\
        0.75 & 550 & 661 & 1122 \\
        1.00 & 550 & 649 & 1052 \\
        \bottomrule
    \end{tabular}
    }
\end{table}
\begin{table}[H]
    \centering
    \caption{Sioux Falls network, deterministic case (Scenario S2): $(TTT, k)$ and $\Delta_{TTT}$ results as compliance $\gamma$ varies for demand level $N_c = 550$. $TTT$ is expressed in minutes.}
    \label{tab:ttt_results_siouxfalls_s2}
    \resizebox{0.8\textwidth}{!}{
    \begin{tabular}{cc ccc ccc}
        \toprule
        & & \multicolumn{3}{c}{\texttt{BS}} & \multicolumn{3}{c}{\texttt{SD}} \\
        \cmidrule(lr){3-5} \cmidrule(lr){6-8}
        $\gamma$ & $N_c$ & $(TTT_{\min}, k_{\min})$ & $TTT_0$ & $\Delta_{TTT}$ (\%) & $(TTT_{\min}, k_{\min})$ & $TTT_0$ & $\Delta_{TTT}$ (\%) \\
        \midrule
        0.25 & 550 & (3096.32, 59) & 3466.19 & 10.67 & (3069.04, 41) & 3466.19 & 11.46 \\
        0.50 & 550 & (3031.04, 66) & 3466.19 & 12.55 & (3002.54, 41) & 3466.19 & 13.37 \\
        0.75 & 550 & (2976.26, 47) & 3466.19 & 14.13 & (2939.83, 46) & 3466.19 & 15.18 \\
        1.00 & 550 & (3029.18, 76) & 3466.19 & 12.61 & (2935.20, 66) & 3466.19 & 15.31 \\
        \bottomrule
    \end{tabular}
    }
\end{table}
\begin{table}[H]
    \centering
    \caption{Sioux Falls network, deterministic case (Scenario S2): computational times (sec) as compliance $\gamma$ varies for demand level $N_{c} = 550$}
    \label{tab:execution_times_s2}
    \resizebox{0.25\textwidth}{!}{
    \begin{tabular}{cccc}
        \toprule
        $\gamma$ & $N_{c}$ & \texttt{BS}  & \texttt{SD} \\
        \midrule
        0.25 & 550 & 830  & 2543 \\
        0.50 & 550 & 434  & 2070\\
        0.75 & 550 & 1185 & 2000 \\
        1.00 & 550 & 662  & 1958 \\
        \bottomrule
    \end{tabular}
    }
\end{table}
It turns out that in general \texttt{SD} provides better solutions with a lower $TTT_{min}$, but at a higher computational cost, in particular up to twice as much for S1 ($\gamma = 0.5$) and up to 3 times as much for S2 ($\gamma = 0.25$). Fig.\ \ref{fig:opt_conf_determ} compare the traffic conditions obtained without RSUs and with the optimal RSUs placement for scenarios S1 and S2, respectively. For scenario S1, the optimal configuration is obtained using the \texttt{SD} method with $\gamma = 1.00$, corresponding to the deployment of 21 RSUs distributed as shown in Fig.\ \ref{subfig:S1_opt_conf}. For scenario S2, the optimal solution is also obtained with the \texttt{SD} method and $\gamma = 1.00$, requiring 66 RSUs distributed as illustrated in Fig.\ \ref{subfig:S2_opt_conf}.
For the four cases shown in Fig.\ \ref{fig:opt_conf_determ}, the total fuel consumption was computed using \eqref{eq:fuel_cons}, \eqref{eq:VSP} and \eqref{eq:total_fuel_cons}, from the beginning of the simulation until the network was completely cleared. The obtained results are reported in Table \ref{tab:fuel_cons}. The results indicate that, under full driver compliance ($\gamma = 1.00$), the optimized RSUs placement reduces the total fuel consumption by approximately 14 L in scenario S1 and 9 L in scenario S2 with respect to the corresponding case without RSUs. These reductions are a direct consequence of the lower congestion levels and shorter travel times achieved by the optimized traffic management strategy.
Besides improving traffic efficiency, the reduction in fuel consumption also implies a corresponding decrease in vehicle exhaust emissions, highlighting the potential environmental benefits of the proposed RSUs deployment strategy.

\begin{table}[H]
   \centering
   \caption{Sioux Falls network: total fuel consumption for scenarios S1 and S2.}
   \label{tab:fuel_cons}
   \resizebox{0.5\textwidth}{!}{
   \begin{tabular}{lcc}
   \toprule
   Scenario & RSUs placement & $J_F$ [L]  \\
   \midrule
   S1 & no RSUs & 92.72  \\
   S1 & optimal (21 RSUs) & 78.94 (-14.86\%)  \\
   \midrule
   S2 & no RSUs & 95.05  \\
   S2 & optimal (66 RSUs) & 86.18 (-9.33\%)  \\
   \bottomrule
   \end{tabular}
   }
\end{table}

\subsubsection{Sioux Falls: stochastic demand case}
\label{sec:siouxfalls_stoc}
While the deterministic scenarios provide a convenient benchmark for evaluating the proposed optimization algorithm, they do not capture the inherent variability of traffic demand. To assess the robustness of the proposed methodology under more realistic conditions, we therefore introduce a stochastic demand scenario.

As in the deterministic case, traffic demand is generated using the residential and work area classification derived from \cite{chakirov2014enriched}. Instead of assigning fixed OD pairs, we associate each network junction with an origin weight and a destination weight according to its functional role within the network. Specifically, residential junctions are assigned an origin weight of 3 and a destination weight of 0.5, whereas work junctions are assigned an origin weight of 0.5 and a destination weight of 3, consistently with the morning home-to-work commuting period. Mixed junctions are assigned unit weights for both origins and destinations. The resulting origin and destination weights are normalized separately over all network junctions to obtain the corresponding origin and destination probability distributions. These distributions are then used to independently sample the origin and destination of each vehicle while ensuring that the two junctions do not coincide. As a result, each simulation generates a different traffic demand realization. The deterministic scenario S1 can therefore be regarded as a particular instance of this more general stochastic framework.

Unlike the deterministic tests, where a single simulation is sufficient to evaluate a candidate RSU configuration, different realizations of the stochastic demand produce different values of the $TTT$. Therefore, each RSU configuration is evaluated over 50 independent simulation runs, and the optimization algorithm uses the average $TTT$ over all runs as the objective function.

For the stochastic case, unlike the deterministic one, we set the compliance level $\gamma$ to 25\%. This choice reflects the realistic assumption that most vehicles are not equipped with I2V communication technology and that, even when available, human drivers generally do not follow routing advices provided by infrastructure, but rather act according to their own interests. 

Table \ref{tab:ttt_results_siouxfalls_s4} presents the $TTT_{\min}$ results for the stochastic case using both the \texttt{BS} and \texttt{SD} methods, whereas the corresponding computational times are shown in Table \ref{tab:execution_times_s4}.

\begin{table}[htbp]
    \centering
    \caption{Sioux Falls network, stochastic case: $(TTT, k)$ and $\Delta_{TTT}$ results for compliance value $\gamma = 0.25$ and demand level $N_{c} = 550$. $TTT$ is expressed in minutes.}
    \label{tab:ttt_results_siouxfalls_s4}
    \resizebox{0.85\textwidth}{!}{
    \begin{tabular}{cc ccc ccc}
        \toprule
        & & \multicolumn{3}{c}{\texttt{BS}} & \multicolumn{3}{c}{\texttt{SD}} \\
        \cmidrule(lr){3-5} \cmidrule(lr){6-8}
        $\gamma$ & $N_c$ & $(TTT_{\min}, k_{\min})$ & $TTT_0$ & $\Delta_{TTT}$ (\%) & $(TTT_{\min}, k_{\min})$ & $TTT_0$ & $\Delta_{TTT}$ (\%) \\
        \midrule
        0.25 & 550 & (1949.21, 57) & 2138.67 & 8.86 & (1947.72, 71) & 2138.67 & 8.93 \\
        \bottomrule
    \end{tabular}
    }
\end{table}
\begin{table}[H]
    \centering
    \caption{Sioux Falls network, stochastic case: computational times (sec) for compliance value $\gamma = 0.25$ and demand level $N_{c} = 550$}
    \label{tab:execution_times_s4}
    \resizebox{0.25\textwidth}{!}{
    \begin{tabular}{cccc}
        \toprule
        $\gamma$ & $N_{c}$ & \texttt{BS}  & \texttt{SD} \\
        \midrule
        0.25 & 550 & 5881 & 34464 \\
        \bottomrule
    \end{tabular}
    }
\end{table}
We note that the best $TTT_{\min}$ is provided by \texttt{SD}, even though it requires many more RSUs (71 compared to 57 for \texttt{BS}). Furthermore, \texttt{SD} takes almost 6 times longer than \texttt{BS} to obtain a solution that is only about 2 seconds better.

\section{Conclusions} \label{sec:conclusion}

This paper addressed the problem of optimal RSUs placement with the objective of minimizing the $TTT$. We proposed a novel optimization framework that integrates a second-order microscopic traffic model with a discrete location ILS metaheuristic in a closed-loop architecture. The traffic model is based on a dynamic follow-the-leader formulation that captures realistic vehicle interactions. RSUs cooperatively collect speed data from near vehicles and return travel time information. The RSUs placement problem is solved with a nested optimization scheme where an outer module determines the optimal number $k$ of RSUs via either bisection search (\texttt{BS}) or stepwise decrement (\texttt{SD}), while the inner module optimizes the location on the road links given each $k$. Two ranking criteria for links are used, path road coverage and link OD coverage that capture both the topological importance of road segments and their potential as traffic diversion points. The framework also integrates an instantaneous fuel consumption model, enabling the assessment of the environmental impact of each RSUs configuration.

The methodology was validated on two networks. On the small-scale network (i.e., Diamond network), we were able to verify the convergence of the proposed algorithms to the same optimal $TTT_{\min}$, achieving up to 36\% $TTT$ reduction under full compliance, with \texttt{BS} proved the be the most efficient strategy, running $5\times$ to $13\times$ faster than \texttt{ES}. On the larger network (i.e., Sioux Falls network), the best deterministic configurations achieved $TTT$ reductions of approximately 20\% (S1) and 15\% (S2), with fuel consumption decreasing by 14.9\% and 9.3\%, respectively. In the stochastic demand scenario with only 25\% driver compliance, the methodology attained an 8.9\% $TTT$ reduction, confirming robustness under realistic conditions. Several insight emerges from the test results. First, RSUs effectiveness strongly depends on the compliance parameter, with the the gap between best and worst configurations decreasing as $\gamma$ decreases. Second, $TTT$ improvement decreases as the network approaches saturation. Third, RSUs location is a highly nonlinear problem, and even augmenting the number of RSUs can lead to worsening of performances. Moreover, locating RSUs on the highest-traffic roads can be counterproductive, as diverting vehicles away from monitored segments may degrade future travel-time estimates. This trade-off between information collection and traffic redistribution motivates the dual-criteria ranking at the core of the proposed ILS.

Future work includes extending the framework to multi-class traffic, exploring dynamic RSU deployment, and incorporating sustainability metrics inside the optimization procedure.
\appendix 
\section{Pseudocodes of the optimization module.}
\label{sec:appendix}
\begin{algorithm}[H]
\scriptsize
\DontPrintSemicolon

\KwData{$G$ (network graph), $M_{od}$ (origin-destination matrix), $k$ (number of RSUs), $\tau_{\max}$ (maximum non-improving iterations)}
\KwResult{Best evaluated RSU configuration $\mathbf{x}_k^*$ for a given $k$ and its corresponding $\delta_{TTT}(\mathbf{x_k^*})$}

Initialize $\mathbf{x} \leftarrow \texttt{GetInitPos}(G, M_{od},  k)$ \textit{// \hyperref[par:init]{Section \ref{par:init}}}\;
$\mathcal{V} \leftarrow \emptyset$, $\delta_{TTT}(x_k^*) \leftarrow - \infty$, $c_{\text{fail}} \leftarrow 0$, $n_k \leftarrow 0$\;

\While{$c_{\text{fail}} < \tau_{\max}$}{
    $(TTT_{\mathbf{x}},  N_{\text{cars}}, P) \leftarrow \texttt{SimulateTraffic}(G, M_{od}, \mathbf{x})$ \textit{ // \hyperref[sec:traffic_simulation]{Section \ref{sec:traffic_simulation}}}\;

    $\delta_{TTT}(\mathbf{x}) \leftarrow \frac{TTT_0 - TTT_{\mathbf{x}}}{TTT_0}$\;

    \If{$\delta_{TTT}(\mathbf{x}) > \delta_{TTT}(\mathbf{x_k^*})$}{
        $\mathbf{x}^*_k \leftarrow \mathbf{x}$\;
    }

    $L ^r \leftarrow \texttt{RankLinks}(G, M_{od}, N_{\text{cars}}, P, \alpha)$ \textit{// \hyperref[par:road_seg_rank]{Section \ref{par:road_seg_rank}}}\;
    $\mathbf{x}' \leftarrow \texttt{SelectTopK}(L^r, k)$ \textit{// \hyperref[par:road_seg_rank]{Section \ref{par:road_seg_rank}}}\;

    \If{$\mathbf{x}' \in \mathcal{V}$}{
        $\mathbf{x}_{\text{pert}} \leftarrow \texttt{Perturbation}(\mathbf{x'}, m, \mathcal{V})$ \textit{// \hyperref[alg:perturbation]{Algorithm \ref{alg:perturbation}}} \;
        $\mathbf{x} \leftarrow \texttt{LocalSearch}(\mathbf{x}_{\text{pert}},  m, \mathcal{V})$ \textit{ // \hyperref[alg:localSearch]{Algorithm \ref{alg:localSearch}}}\;
        $c_{\text{fail}} \leftarrow c_{\text{fail}} + 1$\;
    }
    \Else{
        $\mathbf{x} \leftarrow \mathbf{x}'$\;
    }

    $\mathcal{V} \leftarrow \mathcal{V} \cup \{\mathbf{x}\}$\;
}

\Return $(\mathbf{x}^*_k, \delta_{TTT}(\mathbf{x_k^*}))$

\caption{FindBestConfig(k)}
\label{alg:rsu_optimizer}
\end{algorithm}
\begin{algorithm}[H]
\scriptsize
\DontPrintSemicolon
\caption{Perturbation}
\label{alg:perturbation}

\KwData{$L$, $\mathbf{x}$, $m$, $\mathcal{V}$}
\KwResult{$\mathbf{x}^{pert}$, $R_x^{mobile}$, $R_x^{fixed}$}

Extract indices of active RSUs:
$R_{\mathbf{x}} \leftarrow \{ i \in L : x_i = 1 \}$ \;

\Repeat{$\mathbf{x}^{pert} \notin \mathcal{V}$}{

    Partition $R_{\mathbf{x}}$ into 
    $\{R_{\mathbf{x}}^{fixed}, R_{\mathbf{x}}^{mobile}\}$ 
    such that $|R_{\mathbf{x}}^{mobile}| = m$ \;

    Identify available segments
    $L_{\mathbf{x}}^{empty} \leftarrow \{ i \in L : x_i = 0 \}$ \;

    Select a random subset 
    $R_{new} \subseteq L_{\mathbf{x}}^{empty}$ 
    such that $|R_{new}| = m$ \;

    Construct candidate $\mathbf{x}^{pert}$ such that
    $R_{\mathbf{x}}^{pert} =
    R_{\mathbf{x}}^{fixed} \cup R_{new}$ \;
}

\Return $\mathbf{x}^{pert}, R_{\mathbf{x}}^{mobile}, R_{\mathbf{x}}^{fixed}$ \;

\end{algorithm}

\begin{algorithm}[H]
\scriptsize
\DontPrintSemicolon
\caption{Local Search}
\label{alg:localSearch}

\KwData{$\mathbf{x}^{pert}$, $R_{\mathbf{x}}^{mobile}$,$R_{\mathbf{x}}^{fixed}$, $\mathcal{V}$, $TTT_0$, $S_{max}$}
\KwResult{Local best evaluated configuration $\mathbf{\bar{x}}$}

Initialize neighborhood set:
$\mathcal{N} \leftarrow \emptyset$ \;
Initialize current solution:
$\mathbf{x} \leftarrow \mathbf{x}^{pert}$ \;
\ForEach{$u \in R_{\mathbf{x}}^{mobile}$}{

    Identify available adjacent links:
    $\mathcal{C}_u \leftarrow
    \{ j \in Adj(u) : x_{pert,j}=0 \}$ \;
}

Generate all feasible combinations:
$\mathcal{S} \leftarrow
\mathcal{C}_1 \times \cdots \times \mathcal{C}_m$ \;

\ForEach{$(j_1,\dots,j_m) \in \mathcal{S}$}{

    \If{$|\{j_1,\dots,j_m\}| < m$}{
        continue \;
    }

    Construct candidate configuration
    $\mathbf{x}_{cand}$ using
    $R_{\mathbf{x}}^{fixed} \cup \{j_1,\dots,j_m\}$ \;

    \If{$\mathbf{x}_{cand} \notin \mathcal{V}$
    and
    $\mathbf{x}_{cand} \notin \mathcal{N}$}{

        $\mathcal{N} \leftarrow
        \mathcal{N} \cup \{\mathbf{x}_{cand}\}$ \;
        $\mathcal{V} \leftarrow
        \mathcal{V} \cup \{\mathbf{x}_{cand}\}$ \;

        \If{$|\mathcal{N}| \ge S_{max}$}{
            break \;
        }
    }
}

Evaluate each $\mathbf{x} \in \mathcal{N}$ by computing $\delta_{TTT}(\mathbf{x})$

\Return
$\mathbf{\bar{x}} =
\arg\max_{\mathbf{x}\in\mathcal{N}}
\delta_{TTT}(\mathbf{x})$ \;

\end{algorithm}
\begin{algorithm}[H]
\scriptsize
\DontPrintSemicolon
\KwData{$G$ (network graph), $M_{od}$ (origin-destination matrix), $k_{\min}$ (minimum units), $k_{\max}$ (maximum units), $I_{\max}$ (maximum iterations), $\tau_{\max}$ (maximum non-improving iterations), $\kappa$ (integer decrement step)}
\KwResult{Global best evaluated RSU configuration $\mathbf{x}^*$, $TTT_{\min}$, $\Delta_{TTT}$} 
Initialize $k \leftarrow k_{\max}$, $\Delta_{TTT} \leftarrow - \infty$ \;
$TTT_0 \leftarrow \texttt{SimulateTraffic}(G,  M_{od}, \mathbf{x}_0)$ \textit{ // \hyperref[sec:traffic_simulation]{Section \ref{sec:traffic_simulation}}}\;

\While{$i < I_{\max}$}{
    $(\mathbf{x}_k^*,\delta_{TTT}(\mathbf{x_k^*}))
    \leftarrow
    \texttt{FindBestConfig}(G,k,\tau_{\max})$  \textit{ // \hyperref[alg:rsu_optimizer]{Algorithm \ref{alg:rsu_optimizer}}}\;
    $i \leftarrow i + \tau_{\max}$ \; 
    \If{$\delta_{TTT}(\mathbf{x_k^*}) > \Delta_{TTT}$}{
        $\mathbf{x}^* \leftarrow \mathbf{x}_k^*$\;
        $TTT_{\min} \leftarrow TTT_{\mathbf{x^*}}$ \;
        $\Delta_{TTT} \leftarrow \delta_{TTT}(\mathbf{x^*})$ 
    }

    $k \leftarrow \max(k- \kappa,k_{\min})$\;
}

\Return $\mathbf{x}^*$, $TTT_{\min}$, $\Delta_{TTT}$ 

\caption{Stepwise Decrement (\texttt{SD}) Strategy}
\label{alg:stepwise}
\end{algorithm}
\begin{algorithm}[H]
\scriptsize
\DontPrintSemicolon

\KwData{$G$ (network graph), $M_{od}$ (origin-destination matrix), $k_{\min}$ (minimum units), $k_{\max}$ (maximum units), $I_{\max}$ (maximum iterations), $\tau_{\max}$ (maximum non-improving iterations)}
\KwResult{Global best evaluated RSU configuration $\mathbf{x}^*$,  $TTT_{\min}$,  $\Delta_{TTT}$}
$\mathcal{E} \leftarrow \emptyset$ \textit{// Set of evaluated configurations} \;
Initialize $\Delta_{TTT} \leftarrow - \infty$, $i \leftarrow 0$\;
$TTT_0 \leftarrow \texttt{SimulateTraffic}(G, M_{od}, \mathbf{x}_0)$ \textit{ // \hyperref[sec:traffic_simulation]{Section \ref{sec:traffic_simulation}}}\;
\BlankLine
\ForEach{$k \in \{k_{\min}, k_{\max}, \lfloor(k_{\min} + k_{\max})/2\rfloor\}$}{
    $(\mathbf{x}_k^*, \delta_{TTT}(\mathbf{x_k^*})) \leftarrow \texttt{FindBestConfig}(G, k, \tau_{\max})$ \textit{ // \hyperref[alg:rsu_optimizer]{Algorithm \ref{alg:rsu_optimizer}}}\;
    $\mathcal{E} \leftarrow \mathcal{E} \cup \{(k, \mathbf{x}_k^*, \delta_{TTT}(\mathbf{x_k^*}))\}$\;
    $i \leftarrow i + \tau_{\max}$ \;
    \If{$\delta_{TTT}(\mathbf{x_k^*}) > \Delta_{TTT}$}{
        $\mathbf{x}^* \leftarrow \mathbf{x}_k^*$\;
        $TTT_{\min} \leftarrow TTT_{\mathbf{x^*}}$  \;
        $\Delta_{TTT} \leftarrow \delta_{TTT}(\mathbf{x^*})$ 
    }
}

\BlankLine
\While{$(k_{\max} - k_{\min}) > 1$ \textbf{and} $i < I_{\max}$}{
    $(k_1, k_2) \leftarrow \texttt{BestTwo}(\mathcal{E})$ \textit{// Select top two $k$ based on best $\Delta_{TTT}$}\;
    $k_{\text{mid}} \leftarrow \lfloor(k_1 + k_2)/2\rfloor$\;

    \If{$k_{\text{mid}} \notin \mathcal{E}$}{
            $k \leftarrow k_{\text{mid}}$ \;
            $(\mathbf{x}_k^*, \delta_{TTT}(\mathbf{x_k^*})) \leftarrow \texttt{FindBestConfig}(G, k, \tau_{\max})$  \textit{ // \hyperref[alg:rsu_optimizer]{Algorithm \ref{alg:rsu_optimizer}}}\;
            $\mathcal{E} \leftarrow \mathcal{E} \cup \{(k, \mathbf{x}_k^*, \delta_{TTT}(\mathbf{x_k^*}))\}$\;
            \If{$\delta_{TTT}(\mathbf{x_k^*}) > \Delta_{TTT}$}{
                $\mathbf{x}^* \leftarrow \mathbf{x}_k^*$\;
                $TTT_{\min} \leftarrow TTT_{\mathbf{x^*}}$ \;
                 $\Delta_{TTT} \leftarrow \delta_{TTT}(\mathbf{x^*})$ 

        }
    }
    
    \BlankLine
    $k_{\min} \leftarrow \min(k_1, k_2)$\;
    $k_{\max} \leftarrow \max(k_1, k_2)$\;
    $i \leftarrow i + \tau_{\max}$ \;
}

\BlankLine
\Return $\mathbf{x}^*$, $TTT_{\min}$,  $\Delta_{TTT}$

\caption{Bisection Search (\texttt{BS}) Strategy}
\label{alg:bisection}
\end{algorithm}
\begin{algorithm}[H]
\scriptsize
\DontPrintSemicolon

\KwData{$G$ (network graph), $M_{od}$ (origin-destination matrix), $k_{\min}$ (minimum units), $k_{\max}$ (maximum units)}
\KwResult{Global best RSU configuration $\mathbf{x}^*$, $TTT_{\min}$,  $\Delta_{TTT}$}
$\mathcal{E} \leftarrow \emptyset$ \textit{// Set of evaluated configurations} \;
Initialize $\Delta_{TTT} \leftarrow - \infty$ \;
$TTT_0 \leftarrow \texttt{SimulateTraffic}(G, M_{od}, \mathbf{x}_0)$ \textit{ // \hyperref[sec:traffic_simulation]{Section \ref{sec:traffic_simulation}}}\;
\BlankLine

\For{$k \leftarrow k_{\min}$ \KwTo $k_{\max}$}{
    \BlankLine
    \ForEach{$\mathbf{x} \in \{ \mathbf{x} \in \{0, 1\}^N \mid \sum_{i=1}^N x_i = k \}$}{
        \BlankLine
        $(\mathbf{x}, \delta_{TTT}(\mathbf{x}))  \leftarrow \texttt{SimulateTraffic}(G, M_{od}, \mathbf{x})$ \textit{ // \hyperref[sec:traffic_simulation]{Section \ref{sec:traffic_simulation}}}\;
        \BlankLine
        \If{$\delta_{TTT}(\mathbf{x}) > \Delta_{TTT}$}{
            $\mathbf{x}^* \leftarrow \mathbf{x}$\;
            $TTT_{\min} \leftarrow TTT_{\mathbf{x^*}}$ \;
            $\Delta_{TTT} \leftarrow \delta_{TTT}(\mathbf{x^*})$\;
        }
    }
}

\BlankLine
\Return $\mathbf{x}^*$, $TTT_{\min}$,  $\Delta_{TTT}$

\caption{Exhaustive Search (\texttt{ES}) Strategy}
\label{alg:exhaustive}
\end{algorithm}
\section*{Acknowledgments}
This study was carried out within the Spoke 7 of the MOST -- Sustainable Mobility National Research Center and received funding from the European Union Next-Generation EU (PIANO NAZIONALE DI RIPRESA E RESILIENZA (PNRR) -- MISSIONE 4 COMPONENTE 2, INVESTIMENTO 1.4 -- D.D. 1033 17/06/2022, CN00000023). This manuscript reflects only the authors’ views and opinions. Neither the European Union nor the European Commission can be considered responsible for them.\\
F. L. I.\ is funded by INdAM--GNCS Project, CUP E53C25002010001, entitled ``Numerical Analysis and Algorithm Design for Coupled Nonlinear Differential Systems in Multiscale Complex Dynamics''.\\
F. L. I.\ is member of the INdAM research group GNCS.

\section*{Declaration of competing interests}
The authors declare that they have no known competing financial interests or personal relationships that could have appeared to influence the work reported in this paper.

\section*{CRediT Author Statement}
\noindent
\textbf{E.C.}: Conceptualization,   Funding acquisition, Project administration, Methodology, Supervision, Writing -- review and editing. \\
\textbf{F.L.I.}:  Conceptualization, Data curation, Investigation, Methodology, Software, Validation, Visualization, Writing -- original draft, Writing -- review and editing.\\
\textbf{A.L.P.}: Conceptualization, Data curation, Investigation, Methodology, Software, Validation, Visualization, Writing -- original draft, Writing -- review and editing.\\
\textbf{G.S.}: Conceptualization,   Funding acquisition, Project administration, Methodology, Resources, Supervision, Writing -- review and editing.
\bibliographystyle{elsarticle-harv}
\bibliography{bibliography}
\end{document}